\documentclass[12pt,a4paper]{article}
\usepackage[margin=1.1in]{geometry}
\usepackage{amsmath,amssymb,amsthm}
\usepackage{mathrsfs}
\usepackage{enumitem}
\usepackage{graphicx}
\usepackage{hyperref}
\usepackage{microtype}
\usepackage{setspace}
\theoremstyle{plain}

\theoremstyle{definition}
\newtheorem{definition}{Definition}[section]
\newtheorem{example}{Example}[section]
\theoremstyle{remark}

\newcommand{\Mod}{\mathrm{Mod}}
\newcommand{\Th}{\mathrm{Th}}

\newcommand{\calA}{{\mathcal A}}
\newcommand{\calB}{{\mathcal B}}
\newcommand{\calC}{{\mathcal C}}
\newcommand{\calL}{{\mathcal L}}

\newcommand{\R}{{\mathbb{R}}}

\newcommand{\N}{{\mathbb{N}}}
\newcommand{\Z}{{\mathbb{Z}}}
\newcommand{\Q}{{\mathbb{Q}}}

\newcommand{\be}{\begin{equation}}
\newcommand{\ee}{\end{equation}}

\usepackage[capitalize]{cleveref}   
\newcommand{\eq}[1]{\begin{align}#1\end{align}}

\usepackage{xcolor}

\title{The Logic of Emergence}

\author{David H. Wolpert \\
 Santa Fe Institute, 1399 Hyde Park Road, Santa Fe, NM, 87501\\
 \emph{and} \\
David C.\ Krakauer\\
Santa Fe Institute, 1399 Hyde Park Road, Santa Fe, NM, 87501
}

\begin{document}

\maketitle

%

\begin{abstract}
``Emergence'' refers to 
the fact that different mathematical languages and theories seem to be required to describe 
scientific theories that operate at different spatial and temporal scales.
Highly circumscribed, special cases of emergence have been investigated using many different formal tools, 
including coarse-graining (over space and / or time), renormalization, dimension reduction, 
broken symmetry, and various formulations of causality. 
Here we extend the mathematics of model theory to unify and extend these earlier investigations.
As a preliminary application of our framework, we 
investigate canonical examples of emergence, including the emergence of
genetics from chemistry, of psychology from neurobiology, and 
of different forms of collective behavior from the laws governing their constituent parts. By exploiting the power 
of model theory we also uncover many aspects of
emergence that were previously obscured. In particular, we uncover some of
the rich hierarchical structure that relates different scientific fields. For example,
we show that scientific fields can be
related via a collection of multiple partial orders, all defined over the same set. On the other hand,
under an alternative relation among scientific fields, one which allows cycles, 
we show that it is possible that the set of all of scientific fields must collapse into a single scientific
field. Going forward, these investigations point to the formal tools that 
might need to be developed in order to arrive at a deeper understanding of this emergent hierarchy of scientific reality.  
\end{abstract}

\newpage

\section{The Problem with Emergence}
\label{sec:problem-emergence}

\subsection{What's been done before}

Anderson's ``More is Different''~\cite{Anderson1972} remains the foundational text for many researchers seeking to investigate the concept of emergence. Anderson's core argument is that the laws governing a system at one scale do not determine the organizational principles at a larger scale, because the larger-scale description involves broken symmetries that are not deducible from the underlying equations. A crystal obeys quantum electrodynamics, but QED does not predict the crystal structure, because the ground state breaks the rotational symmetry of the Hamiltonian. Anderson's account is physical, dynamical, and tied to specific mechanisms (spontaneous symmetry breaking, the renormalization group). It has been productive within condensed matter physics and chemistry but struggles as a general theory of emergence across the sciences whose relationships are extremely complicated,  approximate, and probabilistic.

Taking a more philosophical rather than mathematical approach, Chalmers~\cite{Chalmers2006} introduced a distinction between weak and strong emergence. Weak emergence obtains when the higher-level phenomenon is unexpected given the base-level laws but is in principle deducible from them. Strong emergence obtains when the higher-level phenomenon is not deducible even in principle, and new, ``transphysical'' laws are required. Consciousness was Chalmer's preferred candidate for strong emergence, where the physical facts do not fix the phenomenal or psycholgical facts. Strong emergence presupposes a notion of deducibility that is rarely made precise, and it tends to mix ontology (what exists), epistemology (what can be known), and experimental methods (what can be measured). The idea of property dualism (brain and mind are incommensurate) that underlies strong emergence~\cite{Chalmers1996} is either a complicated variety of Anderson's scaling argument, or  a stronger claim for epistemological irredducability. 

The idea of multiple realizability is the observation that the same higher-level property of a so-called ``special science'' 
(so properties like reproductive fitness, pain, warfare) can be realized by different lower-level substrates~\cite{Putnam1967,Fodor1974}.
It has become a general argument for the autonomy of those special sciences, that they can meaningfully be considered without consideration of the lower-level substrates. This can be phrased in Chalmerian terms as saying that the special sciences are weakly emergent from the lower-level ones, strictly speaking, but they can be considered to be effectively strongly emergent.

This is closely related to the question ``what is a macrostate?'' and has been taken up by several authors seeking to make coarse-graining precise~\cite{chorin2000optimal,beck2009model,Israeli2004,Flack2017}. One approach to this issue grounded in physics partitions a microscopic state space into bins, declares each bin a macrostate, and studies the induced dynamics on bins. Many of these studies have emphasized a commutativity requirement wherein a coarse-graining is adequate when coarse-graining the dynamics gives the same result as first evolving and then 
coarse-graining~\cite{wolpert201710,pfante2014comparison,CarrollParola2024,CarrollSingh2019,pfante2014comparison}. However, many features of various scientific fields have nothing to do with dynamics in the first place. This is especially true of the special sciences.

Beyond coarse-graining, aspects of emergence have been approached from the perspective of criticality and phase transitions
(and therefore renormalization group theory), the behavior of systems in specific limits, particularly of large numbers, and algorithmic or statistical complexity. Jensen's work~\cite{Jensen1998} on self-organized criticality proposed that many complex systems spontaneously tune themselves to a critical point, producing scale-invariant behavior (power laws, long-range correlations) without the fine-tuning of an external control parameter. However, again, it is hard to see how this captures the rich and complicated relationships among the special sciences and their lower-level substrates. Reproductive fitness, pain and computation are far more than simply being phase transitions or critical points of some sort.

Building on work stretching back to Nagel~\cite{Nagel1961}, Palacios~\cite{Palacios2019} has explored the relationship between approximation and idealization in inter-theory reduction, arguing that the philosophical literature on emergence has underestimated the role of limits (thermodynamic limit, continuum limit) in connecting levels.  Ay and collaborators~\cite{Ay2011,AyPolani2008} propose that emergence is measured by the degree to which the joint probability distribution of the system's components cannot be factored into a product of independent marginals, and so the system is more than the sum of its parts in a precise information-theoretic sense. 

\subsection{What is missing in all earlier frameworks on emergence}

It is not immediately obvious how to go about connecting or integrating these several perspectives.  Anderson's account describes symmetry and invariance, Chalmers's deducibility. Multiple realizability is about the many-to-one character of the micro-macro map, whereas  the macrostate literature captures the adequacy of coarse-grainings. Jensen's is concerned with the parsimonious origin of new structures, whereas Ay and Polani seek to quantify the informational irreducibility of wholes to parts. 

Furthermore, there are several core features of emergence that are not touched on in any of those frameworks. Most obviously, the different scientific disciplines have different vocabularies, i.e., the constituent elements of one discipline will in general contain concepts that are wholly absent from another discipline. The foundational tools and concepts of different scientific disciplines can by definition appear incommensurable. 

As an example, classical Newtonian dynamics is all about positions, velocities, and forces, which govern the time courses of systems specified with initial positions and velocities. In contrast, Mendelian genetics is all about genes, alleles, (co-)dominance, recessive traits, phenotypes, etc.  Nothing to do with time, nothing to do with velocities or forces. Nor does classical dynamics have anything to do with phenotypes or alleles. In no sense is one discipline recovered from the other after coarse-graining some state space. They are even concerned with wholly different kinds of phenomena. Newtonian dynamics is all about, well, dynamics, and about systems moving through an infinite sequence of states as time evolves. Yet equilibrium theories of genetics, such as the Hardy Weinberg distribution, don't involve detailed
 trajectories through state space.

Despite these differences, there are features that both Newtonian mechanics and Mendelian genetics do have in common. Crucially, both adhere to the strictures of (first order) logic. The logical AND of two propositions is perfectly well-defined in both of them, for example. Evidently, classical dynamics and Mendelian dynamics are magisteria that overlap -- but only partially. How then can we relate them?

These earlier approaches share a limitation that the framework developed below is designed to either overcome, or at the very least, expose their limitations and highlight what a future research program on emergence should take seriously.  A significant limitation that is immediately apparent is that each presupposes a single shared state space that is merely partitioned, coarse-grained, or tuned to a critical point. They therefore cannot even state three features that are central to the relations among the special sciences. The first is the use of genuinely different vocabularies, and the signature of the emergent discipline need not be a coarsening of the signature of the lower one, as it is not for genetics relative to chemistry. The second is the existence of non-dynamical sciences, such as Mendelian genetics, whose content is not a statement about trajectories but equilibrium configurations. The third is the probabilistic and partial character of many inter-level relations, which requires the convergence of multiple overlapping space and time scales. We need to move beyond a shared microscopic state space.

\subsection{Framing Emergence}
\label{sec:logic}

Some work towards such a macroscopic view of emergence has been attempted earlier, using aspects of model theory~\cite{marker2002model,hodges1997shorter} and other tools of mathematical logic, specifically first-order logic.
Nagel's \emph{The Structure of Science}~\cite{Nagel1961} proposed that inter-theory reduction (a weak form of the emergence relationship between phenomena and theories) consists in deriving the sentences of one theory from those of another, supplemented by ``bridge laws'' connecting the vocabularies. Nagel treats reduction as a logical relation between sets of sentences, not as a relation between materials. The framework is perhaps too restrictive and fails to account for approximation, idealization, and the role of limiting procedures.  Tarski's related work on truth and semantic concepts~\cite{Tarski1944,Tarski1953} established that the meaning of a formal theory is given by its models. This  approach to scientific theories, developed further by Suppes~\cite{suppes1999introduction,Suppes2002} and van Fraassen~\cite{vanFraassen1980}, treats a scientific theory not as Nagel does, as a set of sentences in a formal language (the ``syntactic'' view) but as a class of structures (the ``semantic'' view). Bunge's systematic ontology~\cite{Bunge1977,Bunge1979} attempted a comprehensive formal framework for levels of organization, emergence, and the hierarchy of the sciences, using the language of set theory and lattice theory. Bunge insisted that emergence is an ontological category, not merely an epistemological one, and that it requires a formal account of what it means for properties at one level to be ``novel'' relative to properties at another. 

The framework we develop here takes up the natural scientific theories of emergence and situate them in logical program of Nagel, the semantics in Tarski and Suppes, while  respecting the ontological emphasis in Bunge. The central concept is interpretability which describes the  construction of one formal vocabulary inside another. This concept is general enough to capture certain aspects of coarse-graining
and multiple realizability, as well as Nagelian reduction. 
Here we extend this framework to include new concepts like 
probabilistic emergence, stochastic emergence, approximate emergence, 
multiple emergence (in which one discipline emerges from the combination of multiple other disciplines, not just one), and
observer-relative emergence.
These extensions allow us to formally explore the full breadth of the idea of emergence, in in all of its richness. 
And the framework surfaces many aspects of emergence that have not even been appreciated before, never mind investigated formally.

Central to our approach is connecting many features of emergent phenomena  
with associated formal concepts from model theory. These include:
\begin{itemize}
\item A specific
\textbf{physical system}, which is represented in our framework as a model-theoretic 
``structure''. An example is the structure of $N$ point masses orbiting one another under Newtonian gravity.
Another example is a set of $M$ moving electrons in vacuum that are close to one another.
\item A set of ``laws'', or a \textbf{scientific formalism}. In our
framework, a formalism is represented 
as a ``first-order logic theory''. Any physical system that is governed by a particular formalism is represented
in our framework as an associated structure that satisfies the theory of that formalism.
An example of a scientific formalism is Newton's laws, including the precise functions, relations
and constants of those laws. Another example is classical Hamiltonian dynamics, which involves different
functions, relations and constants from Newton's laws.
\item A \textbf{scientific discipline} comprises the collection of all 
scientific formalisms that are reformulations of one another, a property that we identify in model-theoretic terms as those first-order theories being
``mutually interpretable''. An example 
is the collection of Newton's laws, Lagrangian dynamics, Hamiltonian dynamics, etc., all of which
are formalisms in the scientific discipline of classical dynamics.

Any structure that satisfies one theory $T$ also satisfies any other theory $T'$ mutually interpretable 
with $T$. (This is a simple consequence of the
fact that ``mutual interpretability' is a transitive relation, as described below.) In our approach, 
this is reflected in the fact that 
any physical system that is governed by one scientific formalism in a scientific discipline
is also governed by all the other formalisms in that discipline. 
So in particular, any of those other formalisms can be used to describe that physical
system. 

\item And finally, the emergence of one scientific discipline from another, 
which is formalized in our framework using extensions of the model-theoretic 
concept of non-mutual interpretability. An example is the emergence of classical dynamics from quantum dynamics. 
\end{itemize}

While these maps from aspects of emergence in the sciences to concepts in model theory lays a 
foundation for formalizing emergence, they leave out many other aspects of emergence
that are central to the scientific enterprise.
To capture all the nuances of emergence in the sciences, we 
need to extend basic model theoretic constructions in several
ways. 

One particularly important extension of model theory we introduce is designed to capture propagation of a probability
measure from the state space of a microstate-level description of a system to a probability measure
over the state space of a macrostate-level description of the same system. Loosely speaking, 
in the language of model theory this is
a way to propagate a measure
from the domain of discourse of a structure satisfying a theory $B$ to the domain of discourse of a structure satisfying
a theory $A$ that is interpreted in $B$. This extension is needed for example to formalize how a thermodynamic
description of a physical system emerges from a statistical physics description of the same system.

A closely
related extension of model theory we also need to introduce is stochastic emergence. In this case a probability measure
across a set of physical systems described by a scientific
formalism $B$ leads to a probability distribution over the \textit{laws} in a scientific formalism $A$ that
emerges from $B$. This extension is needed for example to formalize how
a probability distribution over the joint states of the neurons in a brain leads to a 
probability distribution over the behavioral laws of the organism housing that brain.  

Yet another important extension of model theory we need to introduce
arises from the fact that often one scientific formalism emerges from \textit{multiple}
other formalisms, not just one. As a simple example, biomechanics emerges from the combination of
biochemistry and Newtonian mechanics, together with molecular biology. And perhaps less obviously, when an existing feature or trait experiences a change in context, giving rise to an apparently new function: what evolutionary biologists call an exaptation.  Finally, we also
introduce an extension to allow for \textit{approximate} emergence, in which one scientific
discipline doesn't precisely arise from another one -- but does to within some error tolerance.

Because our framework is built from model theory and the associated concept of interpretation, formulation of emergence applies to
familiar examples in the mathematics literature of model-theoretic interpretation, e,g., in our framework rational numbers ``emerge'' from the integers.
Our framework also captures a very narrowly
focused notion of ``emergence'' commonly found in the physics literature, as a commutativity property relating the
dynamics of systems to the dynamics of coarse-grainings of those systems.
Going further, our framework captures how
high-level computer programs ``emerge'' from low-level compiled versions of those
programs verified through denotational semantics, as well as how the ideal gas law is described as emerging
from Newtonian mechanics.
The concept of ``broken symmetry'' from condensed matter physics, invoked in the canonical
paper by Anderson, is
another special case of emergence captured in our framework.

\subsection{Roadmap of this paper}

In the next section we review the definitions of  structures and logical theories from model theory.
These definitions only involve specifications of functions, relations and constants along
with some arbitrary set they're all defined over.  
Arguably this bare-bones set of definitions captures \textit{all} of the mathematics that arise in the sciences; different
scientific fields simply use equations that involve different functions, relations and constants defined
over different sets. 

Then in~\cref{sec:interpretability} we review the standard model-theoretic definitions of interpreting one structure in another,
and of interpreting one logical theory in another. Next, in~\cref{sec:programming}, we work through an 
example of emergence involving computer programming languages. This illustrates some of the challenges
in extending model theory to capture all of the phenomenon of emergence. After that, in~\cref{sec:emergence},
we introduce our core definition of what it means for one scientific discipline to emerge from another, grounded
in the model-theoretic concept of interpretation. In~\cref{sec:extensions} we introduce some of the
extensions we need to add to this core definition of emergence, in order to capture the complexities of emergence
in the sciences, e.g., the extension to allow for stochastic emergence, as described above. 

Having presented our abstract mathematical framework, in~\cref{sec:applications}
we present a detailed exposition of how 
our framework applies to some of the phenomena in the sciences that are collectively characterized as ``emergence''.
Next we make some more high-level comments. In~\cref{sec:complexity} we briefly discuss the computational
complexity of the emergence map. Then in~\cref{sec:engineered-natural} we make some broad comments
about our formulation of emergence and some of the higher-level (``softer'') sciences, like neurobiology,
cognitive science, and social science. We end in~\cref{sec:concluding_typology} with a discussion of 
what our framework reveals about the scientific disciplines and their relation ship with one another.

\section{Model Theory}
\label{sec:model-theory}

Every equation in every textbook or research paper involves functions, relations, and / or constants.
Moreover, either explicitly or implicitly, they all involve the Boolean operators (logical NOT, AND, etc.),
the existential and universal quantifiers, and of course the equals symbol.
The minimal mathematics that involves all such equations is called ``first-order logic'', or
``predicate logic''. Particular instances
of such equations --- specific examples of physical systems that satisfy such equations --- are called
``structures''. Accordingly, the minimal mathematics that captures all equations in
science is first order logic, and the minimal mathematics that captures
all examples of those equations is the theory of structures. The combination of first
order logic and the theory of structures is known as model
theory (a term coined by Tarski).

In this section we present a very cursory overview of those elements of model theory that
we shall exploit in this paper. The reader is directed to the literature for more details, e.g.,~\cite{marker2002model,hodges1997shorter}.


\subsection{First order languages }

While terminology differs throughout the literature,  in first-order logic we usually use the term
\textbf{signature}  to mean a
collection of three different sets of \textbf{extra-logical symbols}, namely a set of function symbols (each with a specified arity), a set of relation symbols (each with a specified arity), and a set of constant symbols. (As an aside, we note that constants are simply zero-arity functions, and an $n$-ary function is a special $(n{+}1)$-ary relation.) 
The extra-logical symbols in a signature carry no intrinsic meaning, and as described below, they are names awaiting interpretation (they are ``signifiers'', in effect). As an example, a signature used to analyze arithmetic might include the symbols $\{0, 1, +, \times\}$. 

%

As mentioned above, every equation in every scientific paper involves some set of extra-logic symbols,
signifying functions, relations and constants. In addition to a set of extra-logical symbols, which might
vary from one scientific paper to another, 
all equations in every scientific paper also involve the same
equality symbol. Moreover, whenever such equations are combined, the standard Boolean logic symbols come into play
(even if only implicitly). To reflect this, we define the \textbf{logical symbols} of first-order logic
as the Boolean functions ($\land, \lor, \lnot$),\footnote{Recall that other common
Boolean functions can be defined in terms of these ones, in particular the material
implication symbol, $\rightarrow$.}the quantifiers ($\forall, \exists$), and the equality symbol ($=$).
It is also always assumed that the set of extra-logical symbols contains an
infinite set of variables, typically written as lower-case latin characters.
Conventionally, an infinite set of left and right parentheses are also included in the set of logical symbols,
not because they are necessary for the formalism, but to help with readability. Similarly, we will use
the symbol `:' to demarcate existential quantifiers at the beginning of a formula from the rest of the formula.

A \textbf{language} $\mathcal{L}$ is an infinite set of strings, given by 
the closure of the operation of concatenating the symbols in a signature and the {logical symbols} 
according to the logical apparatus of first-order logic. In other words, it is the closure under
the recursive rules of an appropriate grammar
for building strings of logical and extra-logical symbols that form \textbf{(well-formed) formulas}. 

As an example, $\exists \exists = \lnot x$ is not a well-formed formula, and so not a string
in any language. On the other hand,
in the signature $\{0, 1, +, \times\}$, $\forall x \exists y (x + y = 0)$ is a well-formed formula,
as is $\forall x \exists y (x + y = y + x)$.
While different languages have different sets of extra-logical symbols, they all share the same set of logical symbols. 

If a variable in a formula is quantified (with either $\forall$ or $\exists$), it a \textbf{bound} variable. If a variable is not bound, then it is \textbf{free}.
If a formula has any free variables then it is \textbf{open}. 
If instead all the variables in a formula are bound, the formula is closed,
and sometimes called a \textbf{sentence}. As an example, in the signature 
$\{0, 1, +\}$, $x + y = 0$ is an open formula, while $\forall x \exists y (x + y = y + x)$ is a sentence.

Here we do not consider second-order model theory, based on second-order logic. While that ultimately might  prove necessary --  and
should be relatively straightforward -- most of the work on model theory contents itself with first order logic. Moreover,
first order logic suffices for our purposes here.

\subsection{Structures}

A  {structure} for a signature $\Sigma$ is where the functions, relations and constant 
of $\Sigma$ are realized as operations on a set. Formally, a \textbf{structure} $\mathcal{A}$  is a \textbf{domain} $A$ 
(a non-empty set) together with actual functions $f^{\mathcal{A}}: A^n \to A$ for each $n$-ary function symbol, actual relations $R^{\mathcal{A}} \subseteq A^n$ for each $n$-ary relation symbol, and actual elements $c^{\mathcal{A}} \in A$ for each constant symbol. It is common to
use the term \textbf{interpretation} to refer
to the map from each extra-logical symbol in a language
to the associated entity in the domain of a structure. (Unfortunately,
this conflicts with the meaning of ``interpretation'' between structures,
introduced below. The context should always make the meaning clear.)

As an example of these definitions, the signature 
$\{+, \times, 0, 1\}$, consisting of two binary functions and two constants, is shared by the 
counting numbers, the integers, the rationals, and any commutative ring. 
There are several ways to distinguish those four mathematical systems, despite the fact that they have the same signature. 
One way is to interpret each of them into a structure. This requires 
specifying domains for each of  them, and identifying the extra-logical symbols in each of them with specific functions, relations and constants in those domains. To illustrate, in all four of those systems the domain is a countable set, and for any element $x$ in the domain, 
$x + 0 = x$. On the other hand, for the case of the integers, for any $x$ in the domain there exists some $y$ such that $x + y = 0$, but this is not true for the case of the counting numbers. 

Note that we can interpret each of those four systems into an infinite set of structures. For example, let
$A$ be an arbitrary countably bi-infinite set, with elements $\{\ldots, -1, 0, 1, \ldots\}$. 
That set can serve as the domain of the integers, where the symbol $+$
is interpreted as addition in the usual way. Alternatively, let $Y$ be the countably bi-infinite set with elements $\{\ldots, e^{-1}, e^0, e^1, \ldots\}$. This set as well could serve as the domain of the integers, so long as we now interpret $+$ as multiplication in the usual way.
This provides just two of the infinite set of structures one can use to interpret the language of the integers.

An important distinction to bear in mind is that between a structure and its domain. For example,
we can define a structure of the real numbers, using signature $\{0, 1, \times, +, (\cdot)^{-1}, -(\cdot)\}$, and
domain $\R$. That domain is distinct from the full structure; it is just a set, without any associated
properties.

An important aspect of scientific theories is that they typically involve multiple kinds of variable, e.g., positions,  
temperatures, pressures, categorical states (e.g., sleeping vs. waking vs. exercising), genotype frequencies and fitness values. 
A convenient way to capture this is with a \textbf{many-sorted} structure and an associated language. In such a structure
the domain is partitioned into different subsets, each identified with a distinct \textbf{sort}. The maps over the domain
corresponding to the extra-logical symbols (i.e., the functions, relations and constants) are modified so that each of their
arguments specify which sort it comes from. Similarly, the extra-logical symbols in a many-sorted language
specify the sorts of the associated arguments. In addition, the infinite set of variables of the language are
replaced by a collection of infinite sets of variables, one set for each sort.
(Many-sorted structures can be reduced to single-sorted ones, by introducing unary predicates as sort labels, but doing this obscures the content.)

To recapitulate, every \textit{equation} in any research paper or textbook involves languages, either explicitly or implicitly. 
In other words, it is a WFF.
On the other hand, 
the \textit{examples} presented in every such research paper or textbook has a domain over which those 
functions, relations and / or constants are defined, either explicitly or implicitly. 

However, all those research papers and / or textbooks need a way to reason about the functions, relations and constants,
and associated formulas. They need a way to be able combine those functions, etc., and investigate the consequences of those
combinations. More generally, they need a way to specify a set of laws that apply to all
such examples / structures. In the common vernacular,, they need a ``scientific theory''. This is where logic comes in, and it provides
the stitching needles and the glue for reaching conclusions about structures and theories. These needles and glue are the topic of the next subsection.

\subsection{Theories and Models}

Consider a formula $\phi(\bar a)$ in a language $\calL_A$, where $\bar a = (a_{i_1},\ldots , a_{i_m})$ 
is shorthand for a set of $m \ge 0 $ free variables.
Let $\calA$ be a structure,
with functions, relations and constants defined for every extra-logical symbol in $\calL_A$. 
In general, if we adopt the standard interpretation of the Boolean logical symbols, then
$\phi(\bar a)$ may or may not hold in $\calA$ for some specific tuple $\bar a$. If it does, we write 
\be
\calA \models \phi(\bar a), 
\ee
and say that $\phi(\bar a)$ is \textbf{true} in structure $\calA$.
Alternatively, we sometimes say that $\calA$ \textbf{satisfies} or \textbf{models} $\phi(\bar a)$. If $\phi$
is a sentence, we simply write ``$\calA \models \phi$''.

Given a set of sentences $\{\varphi_i\}$, we say that a structure satisfies that set, or is a model of that set, if it
satisfies each of those sentences individually. 
We say that an additional sentence $\xi$ is a \textbf{(logical) consequence} of 
$\{\varphi_i\}$ if every structure that simultaneously 
satisfies every sentence in $\{\varphi_i\}$ also satisfies
$\xi$. Equivalently, we say that $\{\varphi_i\}$
\textbf{proves} $\xi$, and write $\{\varphi_i\} \models \xi$.

The definition of logical consequence can be quite
awkward to work with, since one needs to consider all structures with given properties.
As a shortcut, a \textbf{deductive system} is a purely syntactic set
of rules for generating all logical consequences of
a given set of sentences. We use the turnstile symbol, $\vdash$, to indicate that one
formula (or set of formulas) implies another under some deductive system. As an example,
suppose that the symbol $\rightarrow$ means material implication, as in propositional logic. Then
if we have both the sentence $\varphi$ and the sentence
$\varphi \rightarrow \xi$,  the ``natural deductive system''
would generate the sentence $\xi$. There are many different
deductive systems, and we won't need to pick out one in this paper.


An $\calL$-$\mathbf{theory}$ $T$ is a set of sentences in the language $\calL$. (Note that
simply by using the logical AND, any finite set of
finite sentences can be combined into a single, finite sentence.)
Often theories of interest will be those for which a set of sentences can serve as a sensible set of \textbf{axioms}, 
from which other sentences can be derived. As an example,
in a language with the binary function symbol $+$ and the constant symbol $0$, a 
possible sentence to include in an associated theory is $\forall x (x+0 = x)$.
Sometimes in the literature
however, the term ``theory'' instead means the set of all logical consequences of the axioms.
As always, the context should make the meaning clear.

As an example, the sentence
$\forall x \exists y : x + y = 0$ can be used as an axiom of the theory of the integers.
However, it is not one of the axioms of the counting 
numbers -- and in fact is not satisfied
in any structure satisfying the axioms of the counting numbers. This illustrates that ultimately, it 
is not the (shared) signature $\{+, \times, 0, 1\}$
that distinguishes mathematical systems, but the
associated axioms and resultant theories that
are expressed in that signature, the resultant truths expressible in it and the domains over which those truths hold.

A theory is not necessarily associated with any particular structure (rather, simply with a language), but it will become 
``semantically meaningful'' when interpreted in that context.
We say that a structure is a \textbf{model} of a theory $T$, or that $\mathcal{A}$ \textbf{satisfies} $T$, and write
$\mathcal{A} \models T$, if every sentence in $T$ holds in $\mathcal{A}$. 

\begin{example}
\label{ex:groups_example}
The language of group theory is the tuple $\{0, \times, \cdot^{-1}\}$, where $\times$ is a binary function and $\cdot^{-1}$ is a unary function. 
The axioms are
\begin{enumerate}
\item $\forall x : x \times 0 = x$;
\item $\forall x: x \times x^{-1} = 0$;
\item $\forall x, y, z : (x \times y) \times z = x \times (y \times z)$;
\end{enumerate}
The language of Abelian group theory is identical, and the axioms are also the same, with the addition of a commutativity axiom:
\begin{enumerate}
\item [4.]$\forall x, y : x \times y = y \times x$;
\end{enumerate}

To illustrate, addition over the integers mod $10$ is a model of the set of all four axioms of commutative group theory (if we interpret
the binary operation of the group as addition mod $10$). On the other hand, the set of $\pi / 2$ rotations about any of the three orthogonal
axes of $\mathbb{R}^3$ are a model of the group axioms, but not of the full set of axioms of commutative groups. 
A logical consequence of the axioms of group theory is, $\forall x, \exists y : y \times (x \times x) = 0$. A
logical consequence of the axioms of commutative group theory, which does not hold for all groups,
is $\forall x, \forall y : y \times (x \times x) = (x \times x) \times y$.
\end{example}


The class of all models of a theory $T$ is written $\Mod(T)$, and is called an \textbf{elementary class}. 
Conversely, the \textbf{complete} theory of a structure, $\Th(\mathcal{A})$, is the set of all first-order sentences true in $\mathcal{A}$. 
Two structures are \textbf{elementarily equivalent} if they have the same complete theory.
We say that a set of sentences is \textbf{consistent} if they have at least one model.\footnote{Note that while common
in model theory~\cite{hodges1997shorter}, this definition is actually weaker than the one used in
first order logic, e.g., in G\"odel's completeness theorem
and his first incompleteness theorems. Those both involve applying a deductive system to generate a set of
sentences $A$ from a set of sentences $B$, and there are no deductive systems in the definition here. However,
by definition,
any deductive system can only derive sentences that are logical consequences of one another. Therefore
the definition in first order logic implies the definition here.} 

A theory $T$ is \textbf{semantically complete} for a particular deductive system, if
for all sentences $\varphi$, $T \models \varphi$ implies 
$T \vdash \varphi$. $T$ is \textbf{sound} for that deductive systems if the converse
holds: $T \vdash \varphi$ implies $T\models \varphi$. 
In contrast, a theory is \textbf{(syntactically) complete} if for every sentence $\varphi$ in its language, 
either $T \models \varphi$ or $T \models \lnot\varphi$. Stated differently, all models of $T$ agree
with one another about the truth of all sentences. 
If a theory is syntactically complete then it is consistent. From now on, unless explicitly
stated otherwise, whenever we refer to a theory as being ``complete'' we mean it 
is syntactically complete.


Most scientific theories are incomplete. For example,
Hamiltonian mechanics admits pairs of physical systems that both satisfy Hamilton's equations, but disagree on other sentences besides those
that Hamilton wrote. (As an illustration, the particles in one system may
be gravitationally bound, while some of the particles
in the other system eventually escape to infinity.)
In general, if we don't explicitly say that a given theory is complete,
then it should be assumed that the theory is a set of axioms from which 
other sentences can be deduced. Syntactic incompleteness is also rife in mathematics.
For example, Euclid's axioms are incomplete, neither establishing the truth or the
falsehood of the parallel line hypothesis.

Let $T$ be any theory that includes the axioms of
Peano arithmetic, and $\varphi$  any sentence that is a logical consequence of those axioms.
G\"odel's \textbf{completeness} theorem for first-order logic guarantees that for the common deductive systems,
$T$ is semantically complete.
However, neither semantic completeness nor soundness means that every sentence $\varphi$ is a logical consequence of 
$T$, i.e., they do not imply syntactic completeness. Indeed, define a set
of sentences $S$ to be ``effectively enumerable'' if there is an associated Turing machine that produces
a list of sentences as it runs, and every one of the sentences in $S$ will eventually appear in that list.
Loosely speaking, G\"odel's
first incompleteness theorem says that $T$ cannot be both consistent and syntactically complete
if it includes Peano arithmetic, and if its axioms are effectively enumerable.\footnote{More precisely,
 G\"odel derived a formulation of this result that only involve a particular deductive system, and which
among other things, allowed for the possibility of syntactic completeness without consistency.
It also only involved ``provability'' rather than truth. Rosser's so-called trick then extends G\"odel's result
to the more commonly stated form of the first incompleteness theorem, which is the one just above.}


\textbf{Isomorphism} between structures is defined in the natural
way, as a bijection between their domains and their signatures
that preserves the interpretations into the domains
of the functions, relations and 
constants. Recall that elementarily equivalent structures are those that satisfy exactly the same sentences. 
In general, if two structures are isomorphic, they
are elementarily equivalent. Moreover, the converse holds if the domain of the structure is finite.
However, it does not hold in general
if the domain of the structure is infinite, due to the
Lowenheim-Skolem theorem. So for example, there
are structures that both satisfy the axioms of Peano arithmetic,
but are not isomorphic. This means there are some structures that model
Peano arithmetic but that differ from the standard structure $\Z$
we all learn in grade school. These ``non-standard'' structures
of arithmetic are very counter-intuitive~\cite{hodges1997shorter,marker2002model}.

For any structure $\mathcal{A}$ with domain $A$, a set $S \subseteq A^n$ is \textbf{definable} if there is an open formula $\varphi(x_1,\ldots,x_n)$ in the language of $\mathcal{A}$ such that $S = \{(b_1,\ldots,b_n) : \mathcal{A} \models \varphi[b_1,\ldots,b_n]\}$. This has exactly the intuitive
meaning for the term ``definable'' one would want: a subset is definable if we can \textit{define} it mathematically using the language $\mathcal{A}$.

To illustrate this concept, return to the axioms of groups, discussed in~\cref{ex:groups_example}.
Suppose we extend the language considered there to include constants $h_1, h_2$. Then in any model of a group,
the left coset of $H := \{h_1, h_2\}$ provides a collection of  one-dimensional definable sets, $\{(x \times h_1) \vee (x \times h_2)\}$.
Similarly, $(x \times y = 0) \vee (x \times y \times y) = 0$ is a two-dimensional definable set. If the group is Abelian, that two-dimensional
definable set is also an equivalence relation.

Note that the $n$ elements of $S$ in the definition of a definable set
are ordered. This allows us to apply the same construction to talk about definable relations, and so 
definable functions. (Functions are definable when the graphs specified by their relations are definable sets.)  Note as well that
whether a set is definable in a structure depends not
only on the domain of that structure, but also on the precise
signature of that structure.
For example, one structure for the real numbers without multiplication
would have signature $(0, 1, +)$, and domain $\R$, the set of reals.
The subset $\R$ given by all elements strictly less than $2$,
$\{x : x < 2\}$, is not definable in that structure. However,
if we expand the signature to $(0, 1, +, \cdot)$, then 
that set is definable, as $\{x : \exists y \; (y \cdot y) + x = 2\}$.

A structure $\calB$ is an \textbf{expansion} of a structure $\calA$ if
their domains are the same, but the signature of $\calA$ is a 
subset of the signature of $\calB$. Alternatively, we say that
$\calA$ is a \textbf{reduct} of $\calB$. We use the terms
``proper expansion'' and ``proper reduct'' accordingly.
If $\calB$ is an {expansion} of  $\calA$, then  the collection 
of definable sets in  $\calB$ is a  superset of the collection of definable sets in  $\calA$. 

Here we extend this concept from structures to theories,
to say that a {theory} $T_B$ is a \textbf{(proper) expansion} of a theory $T_A$
if the set of extra-logical symbols occurring in the sentences
in $T_A$ is a (proper) subset of those occurring in the sentences 
in $T_B$, and in addition, $T_A$ is a (proper) subset of $T_B$. We say that the expansion is consistent if $T_B$ is consistent, and finite if $T_B \setminus T_A$ is finite.



If a theory $T$ is both (syntactically) complete and consistent (and so does not include Peano arithmetic), 
then by definition, all structures that satisfy that theory are elementarily equivalent.
(This does not mean that they are isomorphic in general however, as discussed above.) 
On the other hand, if $T$ is not complete, then in
general there will be pairs of structures that separately satisfy $T$ but are not elementarily equivalent, disagreeing 
with one another about the truth of one or more sentences not contained in either $T$ or the set of logical
consequences of $T$

%

These issues will be important below, when we introduce ``auxiliary constraints'', and then later, when
we introduce ``stochastic emergence''.

\section{Interpretability between structures and between theories}
\label{sec:interpretability}

The central concept of this paper is interpretability. 
Intuitively, interpretability involves a ``translation'', defining
the functions, relations and constants of the first structure (resp., theory) in terms of 
formulas in the second structure (resp., theory). 

We first present the formal definition of this concept for the systematic construction of one structure inside another
that has a different signature. We next extend this to defining the systematic construction of one theory inside another one 
that has a different signature. 


\subsection{Interpretation of structures}

\begin{definition}[Interpretation of structures]
\label{def:struct-interp}
Let $\mathcal{A}$ and $\mathcal{B}$ be structures in signatures $\Sigma_A$ and $\Sigma_B$, having
domains $A$ and $B$, respectively. A \textbf{$k$-dimensional interpretation} of $\mathcal{A}$ in $\mathcal{B}$ consists of:
\begin{enumerate}[label=(\roman*),nosep]
\item A domain formula $\Delta(\bar{x})$ in $\mathcal{L}_B$ with $k$ free variables, picking out $\Delta^{\mathcal{B}} \subseteq B^k$.
\item A formula $E^{\mathcal{B}}(\bar{x},\bar{y})$ in $\mathcal{L}_B$ which 
is an equivalence relation on $\Delta^{\mathcal{B}}$. 
The associated equivalence classes are written as  $\Delta^{\mathcal{B}}/E^{\mathcal{B}}$.
\item For each $n$-ary relation symbol $R \in \Sigma_A$, a formula $\varphi_R$ in $\mathcal{L}_B$, defining a relation on $\Delta^{\mathcal{B}}/E^{\mathcal{B}}$ invariant under $E^{\mathcal{B}}$.
\item For each $n$-ary function symbol $f \in \Sigma_A$, a formula $\varphi_f$ in $\mathcal{L}_B$, defining a function on $\Delta^{\mathcal{B}}/E^{\mathcal{B}}$ invariant under $E^{\mathcal{B}}$.
\item For each constant symbol $c \in \Sigma_A$, a formula $\varphi_c$ in $\mathcal{L}_B$, picking out an $E^{\mathcal{B}}$-class.
\item The induced structure $\mathcal{A}^*_{\mathcal{B}} := (\Delta^{\mathcal{B}}/E^{\mathcal{B}},\, \varphi_R^{\mathcal{B}},\, \varphi_f^{\mathcal{B}},\, \varphi_c^{\mathcal{B}})$ must be isomorphic to $\mathcal{A}$.
\end{enumerate}

\end{definition}

In the sequel, we will typically leave the dimension $k$ implicit. We will sometimes refer to the structure
$\mathcal{B}$ in~\cref{def:struct-interp} as the \textbf{lower} structure, with $\mathcal{A}$ being the \textbf{upper} structure.

\begin{figure}[ht]
\centering
\includegraphics[width=\textwidth]{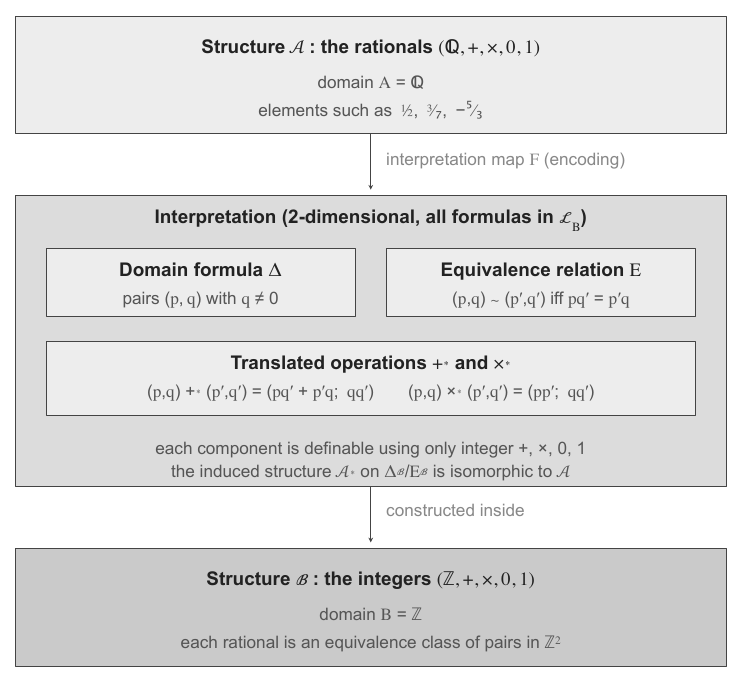}
\caption{The interpretation of the structure $(\mathbb{Q}, +, \times, 0, 1)$ in the structure $(\mathbb{Z}, +, \times, 0, 1)$. The interpretation is 2-dimensional, and each rational maps to an equivalence class of an infinite number of pairs of integers. The domain formula $\Delta$ selects pairs $(p,q)$ with $q \neq 0$. The equivalence relation $E$ identifies pairs encoding the same rational, so that $(1,2) \sim (2,4) \sim (3,6)$. (So each pair in the equivalence class for a rational number $r$ is a numerator and a denominator, respectively, in one particular representation of $r$.) Every translation of an extra-logical symbol in $\mathcal{L}_A$ is a formula in $\mathcal{L}_B$. As an example, the formula in the integers giving the translated version of addition, $+^*$, is defined using only integer addition and multiplication, 
applied to both the numerators and the denominators of the  rationals whose sum it represents. The induced structure $\mathcal{A}^*_{\mathcal{B}}$ on $\Delta^{\mathcal{B}}/E^{\mathcal{B}}$ is isomorphic to $\mathcal{A}$.
}
\label{fig:interpretation}
\end{figure}

The classic example of interpretation of structures
(which has nothing to do with emergence) is how one can use a structure of integers to 
encode a structure of the rationals.
It is illustrative to formulate this example using the concepts of computer programming.
\begin{example}
\label{ex:rationals}
Consider an object-oriented programming language with a single type of variables,  the integers. Suppose the programming
language includes the 
the binary functions $+$ and $\times$ operating on the integers, together with the constants $0$ and $1$ and an infinite set of
parentheses, but no other functions or constants.
Within this programming language, we will define a new object class that we call ``Rational''
This object has two integers as private variables, $a$ and $b$, with $b \neq 0$.
%
Intuitively, we will identify any two objects in this class whose private variable pairs $(a,b)$ and $(c,d)$ are related by
$a\times d = b\times c$ to be equivalent objects, as far as any code using this object class is concerned. To do this,
first adopt the convention that $(a, b)$ is shorthand for an object of class Rational whose two private variables are
$a$ and $b$. Given this shorthand, define
a public method for addition of two objects of the rational class, and a public method for multiplication of 
two objects in class Rational:
\eq{
\label{eq:2}
& PLUS[(a, b), (c, d)] = (a\times d + b\times c,\; b\times d)\\
& TIMES[(a, b), (c, d)] = (a\times c,\; b\times d)
\label{eq:3}
}
In other words, addition and multiplication of Rationals are overloaded binary functions specified in the object class definition,
whose arguments are pairs of integers.

This is an example of interpretation, of the structure of the rational numbers, with signature $\calL_\Q = \{0, 1, \times, +\}$ and
domain $\Q$,
into the structure of the integers, with signature $\calL_\Z = \{0, 1, \times, +\}$ and domain $\Z$, where $k = 2$.\footnote{Note that to agree
with common terminology, in this example we are using the same symbols in each of the two signatures. That has
no formal significance however; at the cost of extra notation we could have written the two signatures 
while distinguishing all of their symbols.}
The equivalence relation $E^{\calB}(\cdot, \cdot)$ from~\cref{def:struct-interp}(ii) is the set of pairs of integers modulo the formula
$a\times d = b\times c$. The interpretation of the two binary functions in the signature of the rationals into
formulas in the signature of the integers is given in~\cref{eq:2,eq:3}. The interpretation of the constants in the
signature of the rationals is simply the associated constants in the signature of the integers.
It is straightforward to verify that these formulas respect the equivalence relation $E^{\calB}(\cdot, \cdot)$, and
that the structure induced by this equivalence relation and those formulas is isomorphic to the structure 
of the rationals.

All of this  is illustrated in Figure \ref{fig:interpretation}. 
An alternative formulation of computer programming in terms of interpretation, involving the compilation
of a program in one language into a program in another language, is presented below in~\cref{sec:programming}. 
\end{example}

\subsection{Interpreting theories}
\label{sec:interpretation_theories}

\begin{figure}[ht]
\centering
\includegraphics[width=\textwidth]{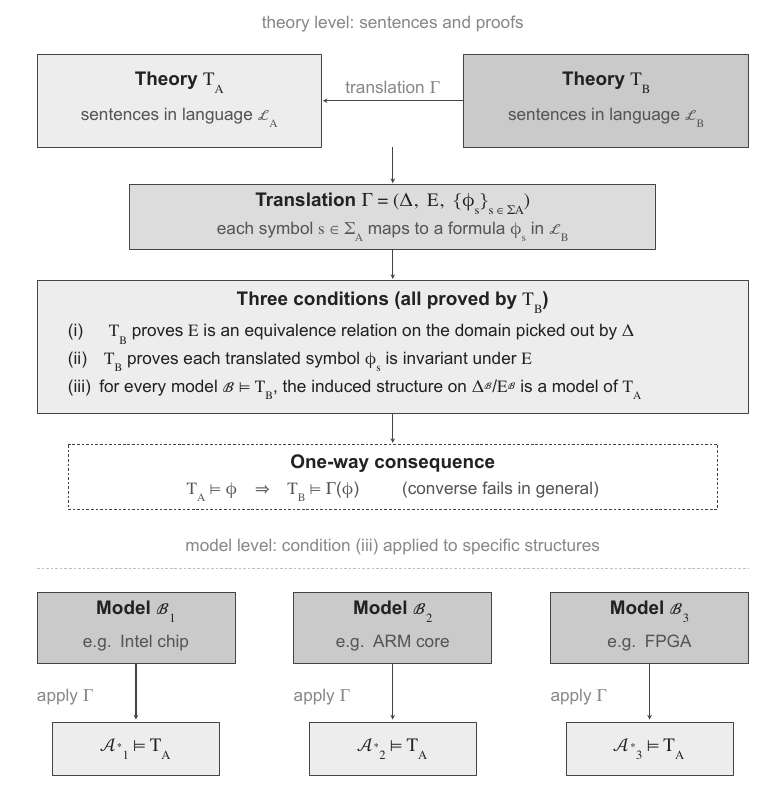}
\caption{Interpretation of theories. The upper half operates at the theory level, where $T_A$ and $T_B$ are sets of sentences in different languages, the translation $\Gamma$ maps each symbol of $\Sigma_A$ to a formula of $\mathcal{L}_B$, and all three well-formedness conditions are proved by $T_B$. The one-way consequence $T_A \models \varphi \Rightarrow T_B \models \Gamma(\varphi)$ holds, and the converse does not. Below the dashed line, condition~(iii) is shown in action, with three different models of $T_B$, each producing via the same $\Gamma$ an induced structure that is a model of $T_A$. The translation works uniformly because it is proved from $T_B$, not verified in any particular model.}
\label{fig:theory-interp}
\end{figure}

The definition of interpreting one structure inside another can be extended to define the interpretation of one theory inside another. The general concept is illustrated graphically
in Figure \ref{fig:theory-interp}, and is formalized as the following pair of definitions:

\begin{definition}(Translations between theories)
\label{def:interp}
A $k$-dimensional \textbf{translation} from 
a theory $T_A$ in language $\mathcal{L}_A$ into a theory $T_B$ in language $\mathcal{L}_B$ is 
a tuple $\Gamma = (\Delta^{{B}}, E^{{B}}, \{\varphi_s\}_{s \in \Sigma_A})$, where 
\begin{enumerate}
\item $\Delta^{{B}}$ is an $\calL_B$-definable relation with
$k$ arguments;
\item $E^{{B}}$ is an $\calL_B$-definable relation with $2k$ arguments;
\item Every $\varphi_s$ is a formula in $\mathcal{L}_B$.    
\end{enumerate}
\end{definition}
\noindent With abuse of notation, we sometimes write $\varphi_s$ as $\Gamma(s)$

Suppose we have a translation 
$\Gamma = (\Delta^{{B}}, E^{{B}}, \{\varphi_s\}_{s \in \Sigma_A})$, 
and that $T_B$ proves that $E^{{B}}$ is an equivalence relation on $\Delta^{{B}}$, i.e.,
on all $k$-tuples that satisfy $\Delta^{{B}}$. Suppose as well that
every $\varphi_s$ is \textbf{invariant} within each equivalence class in $\Delta^{{B}} / E^{{B}}$, i.e. that
for any $\varphi_s$, and every pair of models of $T_B$, $\mathcal{B}$ and  $\mathcal{B}'$,  $\mathcal{B} \models \varphi_s$ iff
$\mathcal{B}' \models \varphi_s$. Then any
model $\mathcal{B}$ of $T_B$ maps $\Delta^{{B}}$ into a definable subset of the domain of that specific model, $\Delta^\mathcal{B}$.
Similarly, it maps $E^{{B}}$ into an equivalence relation over 
$\Delta^\mathcal{B}$, which we write as $E^\mathcal{B}$. Furthermore, 
that model maps each $\varphi_s$ into a formula that is invariant
over the quotient space $\Delta^{{B}}/E^{{B}}$.
So $\mathcal{B}$ induces a structure, whose signature can be identified
with $\calL_A$.

This allows us to make the following definition:
\begin{definition}[Interpretation of theories]
\label{def:theory-interp}
A theory $T_A$ is \textbf{interpretable} in 
a theory $T_B$ if
there is a translation $\Gamma = (\Delta^{{B}}, E^{{B}}, \{\varphi_s\}_{s \in \Sigma_A})$ from $T_A$ into $T_B$ such that:
\begin{enumerate}[label=(\roman*),nosep]
\item $T_B$ proves that $E^{{B}}$ is an equivalence relation on $\Delta^{{B}}$;
\item $T_B$ proves that every $\varphi_s$ is {invariant} within each equivalence class in $\Delta^{{B}} / E^{{B}}$;
\item For every $\mathcal{B} \models T_B$, the associated induced structure on $\Delta^{{B}}/E^{{B}}$ is a model of $T_A$.
\end{enumerate}
We write $T_A \trianglelefteq T_B$ if $T_A$ is interpretable in $T_B$, and refer to $T_B$ as the \textbf{lower} theory, with $T_A$
the \textbf{upper} theory.
\end{definition}
\noindent  In the sequel, we will typically leave the dimension $k$ implicit. 

Interpretability preserves logical consequence in one direction:
\begin{equation}\label{eq:one-way}
T_A \models s \quad \Longrightarrow \quad T_B \models \Gamma(\varphi).
\end{equation}
The converse fails in general. To see this intuitively, note that
strengthening $T_B$ by adding axioms shrinks $\Mod(T_B)$, 
thereby making more sentences provable from $T_B$, without changing what is provable from $T_A$. 
If~\eqref{eq:one-way} were a biconditional, adding axioms to the lower theory could force new consequences in the
upper theory. When the biconditional does hold, the interpretation is called \textbf{faithful}.

Two theories are \textbf{mutually interpretable} if each is interpretable in the other. It is straightforward to confirm that
interpretability is a transitive relation. Therefore mutual interpretability is an equivalence relation.
We also say that two theories are \textbf{bi-interpretable} if their composition is isomorphic to the identity. 

The equivalence classes under mutual interpretability form the \textbf{lattice of interpretability degrees}, 
a dense, distributive partial order that contains incomparable pairs. 
It's important to note though that in his ``linearity conjecture'',
Friedman argues that if we restrict attention
to ``natural'' mathematical theories, then this lattice collapses to a total order~\cite{Friedman2007}.

\subsection{Microstates and interpretation}

It can be useful to be able to measure the ``size'' of the set of elements in
any equivalence class $\Delta^{{B}}/E^{{B}}$.
In the case of interpreting one structure inside another, there is no problem with doing this, and the set of elements
within any particular equivalence class in $\Delta^{{B}}/E^{{B}}$ is perfectly well-defined.
This reflects the fact that the lower structure comes with a pre-defined set of states, its domain.
In contrast, lower theories do {not} come with such a pre-defined domain.
This means that in the case of interpreting one {theory} inside another, there is no 
\textit{a priori} definition of the ``set of elements within any particular equivalence class
in $\Delta^{{B}}/E^{{B}}$''. So we cannot define
the number of elements in such a set, nor define a measure over such a set. 

To circumvent this issue, let $R$ and $R'$ be any two equivalence relations defined
by two formulas in a language $\calL_\calB$, which we write as $\varphi_R$ and $\varphi_{R'}$, respectively.
Assume also that $T_B$ can prove that both of those formulas are 
equivalence relations. Write the associated equivalence classes as the sets $\{r\}$ and $\{r'\}$, respectively.
Then the formula $\varphi_R \land \varphi_{R'}$ is also an equivalence relation, with its equivalence
classes given by the set of all intersections of the two relations' equivalence classes, $\{r \cap  r'\}$.
Similarly, $\varphi_R \vee \varphi_{R'}$ is an equivalence relation, with its equivalence
classes given by the set of all unions of the two relation's equivalence classes, $\{r \cup  r'\}$.
Moreover, every equivalence class defined by $\varphi_R$ is a subset of an equivalence class defined
by $\varphi_R \vee \varphi_{R'}$, and every equivalence class defined by $\varphi_R \land \varphi_{R'}$
is a subset of an equivalence class defined by $\varphi_R$. Therefore we can use subset inclusion
to define a partial order over any set of equivalence classes all of which are both: definable by $T_B$; can be proven by 
$T_B$ to be equivalence classes. 

Suppose we have such a set of equivalence relations, $\mathscr{R}$,
 and that the associated partial order is a finite lattice, with a maximally refined
meet, defined / proven by a finite formula that we write as $\bigwedge_{R \in \mathscr{R}} \varphi_R$. We sometimes refer to the 
equivalence classes of this equivalence relation as \textbf{microstates}. On the other hand, we will refer to the
equivalence classes of
$\bigvee_{R \in \mathscr{R}} \, \varphi_R$ as the \textbf{macrostates} of $T_B$. 
Each of the equivalence classes $\Delta^\calB / E^\calB$ defined in~\cref{def:theory-interp} properly
contains a distinct set of microstates. We will use this set of microstates as the analog, for the case of interpreting 
one theory in another, of the sets of states in each equivalence class that arise when interpreting one structure in another.
The ``size'' of any macrostate is the cardinality of the set of all equivalence classes in that macrostate.

\subsection{General comments on interpretability of theories}

The conditions in~\cref{def:theory-interp} that require $T_B$ to prove the translation $\Gamma$ 
is well-behaved on the quotient $\Delta^B/E^B$ involves
three distinct requirements, each of which must be a logical consequence of $T_B$.

The first requirement is that 
$E^B$ must be an equivalence relation on $\Delta^B$. The formula $E^B(\bar{x},\bar{y})$ defining that equivalence relation
has $2k$ free variables and provides a binary relation on $k$-tuples satisfying $\Delta^{{B}}$. For this relation to 
be an equivalence relation specifically, $T_B$ must prove reflexivity ($\forall \bar{x}\, [\Delta^{{B}}(\bar{x}) \to E^{{B}}(\bar{x},\bar{x})]$), 
symmetry ($\forall \bar{x}\,\bar{y}\, [E^{{B}}(\bar{x},\bar{y}) \to E^{{B}}(\bar{y},\bar{x})]$), and transitivity ($\forall \bar{x}\,\bar{y}\,\bar{z}\, [E^{{B}}(\bar{x},\bar{y}) \land E^{{B}}(\bar{y},\bar{z}) \to E^{{B}}(\bar{x},\bar{z})]$). An arbitrary formula with $2k$ free variables need not define 
an equivalence relation, and $T_B$ must force it to.

The second requirement is that each translated symbol $\varphi_s$ respects $E^{{B}}$. 
If two $k$-tuples are $E^{{B}}$-equivalent,
(they represent the same  macroscopic state), 
then substituting one for the other in any $\varphi_s$ must not change the truth value. Without this, the emergent operations
 would depend on which microscopic representative one happened to examine rather than on the equivalence class, and the macroscopic description
 would leak microscopic information. In the programming example, this says that if two memory configurations encode the same integer, then applying the
 compiled addition routine to either configuration must yield results that are themselves $E^{{B}}$-equivalent.

The third requirement is that the images of the translated function symbols are 
restricted to $\Delta^{{B}}/E^{{B}}$, and in fact for every input class in $\Delta^B / E^B$,
the function must produce an output class that is also in $\Delta^B / E^B$. 
This ensures the induced structure is closed under its own operations.

Because all three of these requirements must be proved by $T_B$, they hold in every model of 
$T_B$, not just in some conveniently chosen one. In the example of~\cref{sec:programming}
below, a compiler that produces correct code on one chip 
but not another has not achieved interpretability. The translation must work across all chips.

If an interpretation from $T_A$ into $T_B$ is non-faithful, and $T_B$ has at least two 
models, then there are models of $T_B$ that disagree on $T_A$, or more
generally (assuming $T_A$ is consistent) disagree on some logical consequences of $T_A$.\footnote{Formally,
let $\calB^1$ and $\calB^2$ be two models of $T_B$, with  $\calB^1_{\Delta^B / E^B}$ and 
$\calB^2_{\Delta^B / E^B}$ the associated induced structures over $\Delta^B / E^B$. If $\Gamma$ is
non-faithful, then there are sentences $\varphi_s$ in $\calL_B$ that both $\calB^1$ and $\calB^2$ satisfy, such
that there are models of $T_A$ that do not satisfy $s$. However, hypothesize that in fact no models of $T_A$
satisfied $s$. Then all such models would satisfy $\neg s$ which by interpretation would mean that
all models of $T_B$ satisfy $\varphi_{\neg s}$, contradicting the fact that both  $\calB^1$ and $\calB^2$
are models of $T_B$ and that they satisfy $\varphi_s$. Therefore the hypothesis is wrong, and
so not all models of $T_A$ agree on $s$.}
Therefore so long as $T_B$ has more than one model, a necessary condition for the interpretation to be non-faithful is that
$T_A$ is not syntactically complete. When $T_A$ is not syntactically complete there are sentences $\phi$ 
in $\calL_A$ that are not logical consequences of $T_A$, and two models
on $\Delta^{{B}}/E^{{B}}$ will disagree on some of those sentences.

\begin{example}
As an illustration of such non-faithfulness, 
suppose that $T_A$ involves two functions $f, g$ that separately obey the same axioms, but no sentences in $T_A$ say how those two functions 
are related. Then the interpretation might result in induced
structures on $\Delta^{{B}}/E^{{B}}$ that don't distinguish 
between those two functions, and therefore disagree 
with one another concerning sentences
in $\calL_A$ that do distinguish those two functions. In this case $T_A$ is not syntactically complete.

As an explicit example, consider a ``two-product group" theory $T_A$, with two binary operations $f, g$ 
and the four axioms:
\eq{
	\exists e \,:\, \forall x, f(x, e) \;=\; x \\
	\exists e \,:\, \forall x, g(x, e) \;=\; x \\
	\forall x \exists y \;:\; f(x, y) \;=\;  e \\
	\forall x \exists y \;:\; g(x, y) \;=\;  e
}
$T_B$ might include an additional axiom specifying that  over
$\Delta^{{B}}/E^{{B}}$ 
that obeys the formula
\eq{
	 \forall x, \forall y \, f(x, y) = g(x, y)
}
However, $T_A$ does not include that axiom, and there are both models of $T_A$ that satisfy
it and models that don't. So the interpretation is non-faithful and $T_A$ is incomplete.
\end{example}

Recall that Peano arithmetic is fundamental to almost all scientific
theories, but is not syntactically complete, by G\"odel's first incompleteness theorem. 
However, this does not mean that any theory $T_A$ 
must by syntactically incomplete if it can be interpreted into another theory $T_B$ which includes
the axioms of Peano arithmetic. 
The fact that $T_B$ is incomplete means that there are 
sentences in $\calL_B$ on which two distinct models of $T_B$
disagree. However, those sentence might not be expressible as formulas 
in $\Delta^{{B}}/E^{{B}}$, the projection of the signature of $T_A$ into $T_B$. So their disagreement
may not back-project into an associated disagreement in $T_A$.

These considerations will be central to the some of the extensions of emergence we introduce below.

\section{From Languages and Machine Code to Hardware}
\label{sec:programming}

The denotational semantics of a minimal, non-optimizing
high-level programming language is an instance of model-theoretic interpretability of one theory
into another.\footnote{A full optimizing compiler is not a model-theoretic interpretation. This is because it is
context-sensitive, operates on whole programs, and preserves only the computed function. The
interpretation is the definitional relationship between the high-level operations and their machine-level
definitions.}   Denotational semantics is used to verify that the meaning of  source code and target code are identical. 
The results of compilation thus provide a natural bridge to the study of emergence~\cite{Krakauer2024}. We develop this idea in detail because it grounds the abstractions introduced above
in familiar mechanisms. Furthermore, the 
relationship between a target program and the hardware that executes it exemplifies the distinction between theory and structure. High
level languages compiled down into machine code that configure time-varying hardware-states cover aspects of the ontological and epistemological
factors described by Nagel, Tarski and Bunge as well as scientifically motivated definitions. 

\begin{figure}[ht]
\centering
\includegraphics[width=\textwidth]{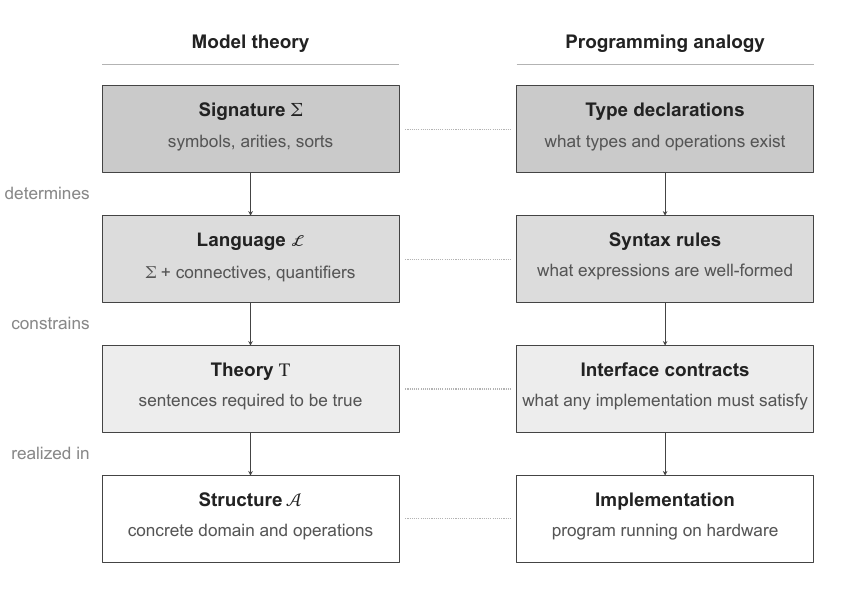}
\caption{The four-level hierarchy of model theory (left) and its programming analogy (right). The signature provides the vocabulary from which the language is built, the theory selects which sentences are asserted, and the structure provides a concrete domain in which the assertions hold. In the programming analogy, type declarations determine what expressions are syntactically valid, interface contracts constrain what any implementation must satisfy, and the implementation is a concrete program that compiles and passes all tests. An interpretation is the construction of a library that implements one interface using only the primitives of another. The emergent theory is the high-level API. And the underlying theory provides the assembly-level operations.}
\label{fig:hierarchy}
\end{figure}

\subsection{Two signatures, two languages}

Consider a minimal high-level language $H$ and a machine-level language $M$. The signature $\Sigma_H$ of $H$ contains a sort for integers, a sort for Booleans, function symbols for arithmetic ($+_H, \times_H$), a conditional function symbol $\mathsf{if}$, and constant symbols $\mathsf{true}, \mathsf{false}, 0_H, 1_H$. The signature $\Sigma_M$ of $M$ contains a sort for bit-strings of length $w$ (words), function symbols for bitwise AND, OR, NOT, addition modulo $2^w$, and a shift operator, and a constant symbol $\mathbf{0}$ (the zero word).

These signatures share no non-logical symbols. The language $\mathcal{L}_H$ can express statements like $\forall x\, (x +_H 0_H = x)$, and the 
language $\mathcal{L}_M$ can express statements like $\forall v\, (\mathsf{AND}(v, \mathbf{0}) = \mathbf{0})$. Neither language can 
directly state sentences in the other's vocabulary.

Let $T_H$ be the high-level theory,  the $\mathcal{L}_H$-sentences derivable from the axioms that define the high-level language, equivalently the properties that hold of every well-formed high-level program. Let $T_M$ be the machine theory, the $\mathcal{L}_M$-sentences derivable from the axioms of the machine. In practice, the first thing a compiler does with a program is check that each of its lines is a syntactically valid sentence of $\mathcal{L}_H$, that is, that the program belongs to the language whose theory is $T_H$.

\subsection{Interpretation as denotational translation}

The denotational translation underlying compilation interprets one theory in another.  This definitional translation provides three components required by Definition~\ref{def:theory-interp}.

The first is a domain formula, and the integers of $H$ are encoded as equivalence classes of bit-strings in $M$. For unsigned integers up to $2^w - 1$ the encoding is one-dimensional, with each integer mapping to one word, while for arbitrary-precision integers it is multi-dimensional, with each integer mapping to a tuple of words representing a linked sequence of digits.

The second is an equivalence relation, and two bit-strings are equivalent when they encode the same integer. For fixed-width arithmetic this is just bitwise equality, while for arbitrary-precision representations it is the relation ``these two linked lists represent the same number.''

The third is a definable translation of each symbol. The high-level $+_H$ is translated into a formula in $\mathcal{L}_M$ that chains together word-level additions with carry propagation. The conditional $\mathsf{if}$ is translated into a formula involving branch-on-zero. Each translation is a formula in $\mathcal{L}_M$, definable from the primitives of $M$, invariant under the equivalence relation.

This translation constructs, inside any model $\mathcal{M}$ of the machine-level theory $T_M$, an induced structure on $\Delta^{\mathcal{M}}/E^{\mathcal{M}}$ that is itself a model of the high-level theory $T_H$. Because this holds for every model of $T_M$ and is provable from $T_M$, we have $T_H \trianglelefteq T_M$.

\subsection{The theory-structure distinction and hardware as model}

The theory $T_M$ specifies how the machine-level operations behave, how addition modulo $2^w$ is associative and commutative, AND distributes over OR, and so on. These are sentences and abstract constraints. and they say nothing about voltage levels, clock speeds, or silicon.

A particular processor (an Intel chip, an ARM core, an FPGA configured to implement the same instruction set) is a structure $\mathcal{M}$ satisfying $T_M$. Different processors are different structures satisfying the same theory. They are elementarily equivalent (they agree on all first-order sentences about word-level operations) but they are not isomorphic, and they differ in physical properties (power consumption, heat dissipation, propagation delay) that are not expressible in $\mathcal{L}_M$.

The same program runs on different hardware. The program is a consequence of the theory and the execution is a process in a structure. The theory constrains the logic and the structure embodies the logic. When we say that the high-level language ``emerges from'' the machine code, we are making a claim at the level of theories ($T_H \trianglelefteq T_M$), not at the level of structures. Different chips give different structures, but the interpretability relation holds uniformly across all of them, because it is proved from $T_M$, not from any particular model. What holds for all models of a theory versus what holds in a particular model is the difference between emergence and instantiation.

The interpretation map (the denotational translation) is one-directional. Every sentence provable in $T_H$ translates to a 
sentence provable in $T_M$. But different choices for $T_M$ all prove sentences that have no counterpart in $T_H$, 
among them overflow behavior and alignment constraints.
The machine-level theory draws distinctions invisible at the high level. This is what it means to say the lower level is ``richer,'' and the notion can be made fully formal. Fix the translation $\Gamma$ with domain formula $\Delta$ and congruence $E$. The emergent theory $T_H$ can speak only of $E$-invariant sets, and a subset of $\Delta$ is visible to $T_H$ exactly when it is a union of $E$-classes of the form $\Gamma(\psi)$ for some $\mathcal{L}_H$-formula $\psi$. To say that $T_M$ is richer than $T_H$ is to say that some $\mathcal{L}_M$-definable subset of $\Delta$ is not a union of $E$-classes, but splits at least one class. Overflow behavior, alignment, and the like
are exactly such $\mathcal{L}_M$-definable predicates, separating memory layouts that lie inside a single integer-class and so are indistinguishable in $\mathcal{L}_H$. A strictly stronger reading is available through faithfulness, and $\Gamma$ is faithful when $T_H \models \varphi$ if and only if $T_M \models \Gamma(\varphi)$, while $T_M$ is richer when $T_H \trianglelefteq T_M$ holds but no interpretation between them is faithful, so that $T_M$ settles translated sentences that $T_H$ leaves open~\cite{VisserFaith2005}.

Conversely, $T_H$ may contain constructs (garbage collection, dynamic typing, closures) that require elaborate construction from the primitives of $T_M$. The higher level is ``concise'' in the sense that short sentences in $\mathcal{L}_H$ compress long formulas in $\mathcal{L}_M$. 
Furthermore, the translation is computable but may involve exponential blow-up in formula length. We return to this issue in Section~\ref{sec:complexity}.


\section{Emergence}
\label{sec:emergence}

\subsection{Emergence as constrained interpretability}

An important feature of the compiler example is that $T_H$ is interpretable in $T_M$ outright, and no additional axioms are needed. In scientific practice the situation is different. 

The ideal gas law is not interpretable in Newtonian mechanics alone since one must add typicality conditions on initial configurations, e.g., restricting to systems with large particle numbers, requiring all collisions be perfectly elastic, defining the shape of the enclosing container, etc. The emergent theory ($PV = NT$) becomes interpretable only in the enriched theory $T_{\text{Newton}} \cup C$, where $C$ encodes these physical constraints as sentences in the lower theory of Newtonian mechanics.

\begin{definition}[Emergence]
\label{def:emergence}
Let $T_A$ and $T_B$ be two theories, and let $\widetilde{T_B} = T_B \cup C$ for a finite, consistent set $C$ of sentences of $\mathcal{L}_B$, such that
$T_A \trianglelefteq \widetilde{T_B}$. We 
then say that
$T_A$ \textbf{emerges from} $T_B$ \textbf{using auxiliary constraints} $C = \widetilde{T_B} \setminus T_B$. 
\end{definition}

Whenever $\widetilde{T_B} \setminus T_B = \varnothing$, emergence reduces to interpretability. More generally, 
adding auxiliary constraints restricts the class of models under consideration. From among all systems obeying $T_B$, 
they select those that exhibit the emergent regularity.

Of course, rather than use~\cref{def:emergence} one could instead simply say that $T_A$ is interpretable in $\widetilde{T_B}$.
The benefit of~\cref{def:emergence}  is that it allows us to capture the case where multiple theories $\{T_i\}$ all emerge from
the same theory $T_B$, but each using different auxiliary constraints, where those auxiliary constraints are not consistent with 
one another. In this case there is no single theory $T$
such all of the  $\{T_i\}$ are interpretable in $T$, even though they all emerge from the same theory. For example,
both the laws of thermodynamics and those of mechanical engineering emerge from Newtonian mechanics, even though neither
is interpretable in Newtonian mechanics.\footnote{This raises an interesting mathematical question.
Suppose we are given a set of theories $\{T_i\}$
and a language $\calB$. Suppose we define the ``size'' of a theory
as the number of sentences in that theory, 
or analogously to the case in proof complexity, as the minimal number of distinct variables  
needed by those sentences. The question is, what is the largest theory $T_B$ in  $\calB$ such that all of the theories in  $\{T_i\}$
emerge from $T_B$? Such a theory would be the ``most macroscopic'' of all theories that
simultaneously are obeyed by all the theories in  $\{T_i\}$.}

If $T_A$ is interpretable in $T_B$ and $T_B$ emerges from $T_C$, then $T_A$ emerges from $T_C$.
So emergence without constraints is transitive.
However, suppose instead that $T_A$ emerges from $T_B$, while $T_B$ is interpretable in $T_C$.
In this case we cannot conclude that $T_A$ emerges from $T_C$. The problem is that the auxiliary constraints added
to $T_B$ (to allow $T_A$ to emerge from $T_B$) might not themselves be interpretable in $T_C$.

Note that there are several different ways that $T_A$ might fail to be interpretable in $T_B$ if we don't add the 
auxiliary constraints. Any of the three parts of~\cref{def:theory-interp} might fail, even though
they hold if we expand $T_B$ by adding the sentences in $\widetilde{T_B} \setminus T_B$. Note as well
that $T_A$ may emerge non-faithfully from $T_B$, if $T_A$ is not
syntactically complete. This is the starting point for the concept
of ``stochastic emergence'' discussed below.

There are numerous ways that we must augment this core definition of emergence in~\cref{def:emergence} to capture all the 
ways the term is used in the sciences and philosophy. For example, classical genetics  as emergent from chemistry, discussed in Section~\ref{sec:genetics-chemistry}, lies at the far end of this spectrum
from the  thermodynamics as emergent from Newton's laws. 
Genetics cannot be made interpretable in chemistry
by adding a set of auxiliary constraints $C$. The reason is simple:
the sentences in any such $C$ are drawn from $\mathcal{L}_{\text{chem}}$, and $\mathcal{L}_{\text{chem}}$
does not contain the extra-logical
symbols that are needed to capture the intrinsically statistical nature of the higher theory of genetics.
The lower theory must itself be augmented with additional mathematical
properties to capture this statistical nature of the higher theory. This is elaborated below in~\cref{sec:emerge_with_measure}
and~\cref{sec:stoch_emergence}.

\subsection{Richness, expressivity, and parsimony}
\label{sec:parsimony}

The preceding definitions provide several 
quantities that are key to the model theoretic framing of emergence, namely the dimension $k$ of the interpretation, 
the size of the equivalence classes under $E$,
and the expansion of each extra-logical symbol in the upper language 
into a string of symbols in the lower language. Together they determine the parsimony of the emergent theory relative to the lower theory.

Consider two structures $\mathcal{A}$ and $\mathcal{B}$ with signatures $\Sigma_A$ and $\Sigma_B$, respectively,
along with domains $A$ and $B$. Suppose that in addition $\mathcal{A}$ emerges from $\mathcal{B}$ using 
a specified $k$-dimensional translation map.      
The dimension $k$ of that translation is the number of $B$-coordinates required to represent a single element of $A$. In the interpretation of $(\mathbb{Q},+,\times,0,1)$ in $(\mathbb{Z},+,\times,0,1)$, we have $k=2$ where each rational requires a pair of integers. This example concerns the interpretation of one structure in another rather than one theory in another. The theory-level version, which requires axiomatizing both the rationals and the integers, is the standard construction carried out in the model-theory literature~\cite{Hodges1997}. 

As another example, in the emergence of the ideal gas law 
from Newtonian mechanics, each macrostate (a quadruple $(N,P,V,T)$) is 
encoded by a configuration of $N$ particles in a $6N$-dimensional phase space, so $k$ is proportional to the number of particles. 
The sentence ``$P = 3$'' in the emergent theory, which names a single macroscopic quantity, translates via $\Gamma$ into a formula 
over $k$-tuples of phase-space coordinates. We say that the dimension $k$ measures the {\bf width} of the translation, and it counts 
how many microscopic degrees of freedom are bundled into one emergent (extra-logical) symbol.

Suppose that the
equivalence class $\bigvee_R \, \varphi_R$ 
discussed at the end of~\cref{sec:interpretation_theories} is well-defined and finite.
So the number of microstates within each equivalence class in $\Delta^{{B}}/E^{{B}}$ is well-defined and finite.
In such a case,
the equivalence classes $[b]_E$ defining $\Delta^{{B}}/E^{{B}}$ provide a different 
measure of compression. Every sentence of $T_A$, when translated via $\Gamma$, is invariant over 
these classes, and it cannot distinguish microstates within a single class.
A large equivalence class means that many $k$-tuples in $B^k$ are identified as representing the same emergent state. 

The number of distinctions the emergent theory can draw is bounded by $|\Delta^{{B}}/E^{{B}}|$, which is smaller than $|\Delta^{{B}}|$ by a factor equal to the average class size. 
Large classes mean the emergent theory is blind to most microscopic detail. 
Assuming some appropriate measure over the number of equivalence classes in 
$|\Delta^{{B}}/E^{{B}}|$, we refer to this expected number of microstates
in each equivalence class as the {\bf depth} of the translation, and write it as $\overline{|[b]_E|}$.
An obvious variant is to instead consider the maximum over all equivalence 
classes of the number of microstates per class, assuming that is well-defined.
Yet another variant arises when probability measures are defined using an extension of $\mathcal{B}$,
as described below in~\cref{def:emergence_measure}. In such a situation the probability
measure allows us to dispense with microstates entirely.

Width and depth characterize features of the translation's map between
the domains of two structures that satisfy $T_A$ and $T_B$, 
respectively. However, the translation
also involves a map between the signatures of those two theories.
Each extra-logical symbol $\alpha \in \Sigma_A$
gets mapped by the translation into a
formula $\varphi_\alpha$, which consists of 
a string of (logical and extra-logical) symbols 
which has length $|\varphi_\alpha|$. So the maximum
over all $\alpha \in \Sigma_A$ of $|\varphi_\alpha|$
is a measure of how much the symbols in the
upper theory compress symbols in the lower theory. 
We call that maximum the \textbf{blow-up} of the translation, and write it as
$\widehat{\{\varphi_{s}\}}$.

Parsimony in the emergent theory $T_A$ relative to $T_B$ is the joint effect of all of these quantities. When $k$ is large, each emergent concept encodes many microscopic degrees of freedom. When the equivalence classes are large, many 
microstates that could be distinguished using the theory $T_B$ are treated as identical. Finally, when 
the blow-up is large, each symbol in $\calL_A$ represents many symbols in $\calL_B$. The product
of these types of compression is
\[
k \;\cdot\; \overline{|[b]_E|} \; \cdot\; \widehat{\{\varphi_{s}\}}
\]
(dimension times average class size times blow-up factor) 
gives a rough measure of the compression ratio, and counts how many microscopic degrees of freedom are collapsed into one emergent degree of freedom. A high compression ratio is what makes the emergent theory so efficient since a short sentence in $\mathcal{L}_A$ encodes, via $\Gamma$, a long and complicated formula in $\mathcal{L}_B$.

This also gives a formal criterion for the asymmetry between ``richer'' and ``expressive'' that the programming example illustrated informally. A theory $T_B$ is richer than $T_A$ when $T_A$ is interpretable in $T_B$ but not vice versa. When there is that asymmetry, there exist $\mathcal{L}_B$-definable sets in models of $T_B$ that are not unions of $E$-classes, and they cut across the classes, distinguishing microstates that the emergent theory cannot. The sets the emergent theory can name are exactly the $E$-invariant ones, namely the $\Gamma$-images of $\mathcal{L}_A$-definable sets, and the distinctions that $T_B$ draws and $T_A$ cannot are precisely those $\mathcal{L}_B$-definable sets that fail to be $E$-invariant. 

The emergent theory $T_A$ is more expressive in the complementary sense, and it can state, in a short sentence using an appropriate vocabulary, a regularity that in $\mathcal{L}_B$ requires quantification over all elements of a large equivalence class. 
Expressivity provides both a compression and a unified mechanism as encoded in signatures and sentences, whereas richness describes the resolution. They are related by the interpretation map as two sides of the same asymmetry.

In the compiler example, the high-level sentence ``$x + y = z$'' compresses a formula in $\mathcal{L}_M$ that involves carry propagation across $w$ bit positions, overflow checking, and so forth. The compression is possible because the equivalence relation $E$ identifies all memory layouts that encode the same integer and the dimension $k$ bundles $w$ bits into one high-level value. The emergent theory (the high-level language) is more concise than the lower theory (machine code) in the sense that it says the same thing in fewer symbols, but the lower theory is richer in the sense that it can distinguish many memory layouts that the high-level theory cannot describe. The same logic of compression and resolution applies to coarse graining, multiple-realizability and symmetry breaking since all of them express the value of effective theories in terms of their advantages beyond detailed descriptions whose partitions are not easily observed.

\subsection{Insufficiency of interpretations with constraints}
\label{sec:genetics-chemistry}

The interpretation of the structure of the
rationals in the integers is a relation between two fixed structures, and its criterion of success is an isomorphism. Interpreting one theory in another is different as per Definition~\ref{def:theory-interp}, avoiding the need for
isomorphism. Indeed, there is no single target structure to be isomorphic to.  The requirement is that every model of the lower theory 
must induce, through one uniform translation, a model of the upper theory. 

Classical genetics in relation to chemistry is an example of this difference between interpretation among structures
and interpretation among theories. It also shows
why plain interpretability is not enough, with or without auxiliary constraints.

Let $T_{\text{gen}}$ be the theory of classical Mendelian genetics, with a many-sorted
signature of a certain number of genes, alleles, loci, genotypes, and the relations of dominance and segregation.
Let $T_{\text{chem}}$ be a theory of chemistry, with a many-sorted signature of atoms, bonds, and molecular configurations,
along with unary functions taking values like hydrophobicity, polarity, charge, etc., and other functions specifying stoichiometrically
allowed reactions, etc.
The two signatures share no non-logical symbols. Moreover, 
a single allele corresponds to an enormous configuration of atoms, so the domain formula $\Delta$ would select large tuples of chemical variables and the congruence $E$ would group those configurations into 
alleles of genes. 

Strict interpretation into $T_{\text{chem}}$ fails immediately. Mendel's law of segregation, under which the two alleles at a locus separate into gametes in equal proportion so that a monohybrid cross yields offspring in a $3\!:\!1$ ratio, is not a theorem of chemistry. No sentence of $T_{\text{chem}}$ proves it, because the ratio is a consequence of the mechanics and combinatorics of chromosome pairing at meiosis, not of the chemistry of any molecule. Thus $T_{\text{gen}} \not\trianglelefteq T_{\text{chem}}$. This is a typical situation where the vocabulary
of an emergent discipline is not a coarse-graining of the vocabulary of the theory it emerges from.

Adding constraints does not repair this. One might hope to adjoin to $T_{\text{chem}}$ a set $C$ of constraints encoding diploidy and the meiotic apparatus and thereby obtain $T_{\text{gen}} \trianglelefteq T_{\text{chem}} \cup C$. 
The obstruction is 
that
definition~\ref{def:theory-interp} requires
$E$ to be $\mathcal{L}_{\text{chem}}$-definable and each symbol of $\Sigma_{\text{gen}}$ to be sent to a single fixed formula of $\mathcal{L}_{\text{chem}}$ that $T_{\text{chem}} \cup C$ proves invariant.
Neither requirement holds here. The equivalence that genetics needs is sameness of allele up to phenotypic effect, and phenotypic effect reaches downstream into a developmental and environmental cascade that lies outside the chemistry of the DNA molecule itself, so $E$ is not definable in $\mathcal{L}_{\text{chem}}$ and cannot be made so by adding
sentences to the set of $\mathcal{L}_{\text{chem}}$-axioms. The characteristic relations are context-sensitive for the same reason. Whether an allele is dominant, and with what expressivity, depends on genetic background, modifier loci, and environment, so there is no fixed, context-free formula $\varphi_{\text{dom}}$ of $\mathcal{L}_{\text{chem}}$ into which the dominance symbol could be translated. This is exactly the compiler's failure. A construct whose meaning depends on its whole context is not sent to a single formula, and a congruence that is not definable in the lower language cannot be recovered by strengthening the lower theory.

The contrast with the compiler is  instructive. There, a stripped-down, non-optimizing high-level language had a core that was interpretable in the same lower theory $T_M$.  Genetics has no such core and nothing in it is interpretable in chemistry alone, whether or not chemistry is constrained. Its interpretable core appears only after the lower theory is expanded.

This expansion is formalized in~\cref{sec:extensions} below
as ``multiple emergence''. Genetics emerges jointly from chemistry, from a theory of chromosomal mechanics and the combinatorics of segregation, which is where the $3\!:\!1$ ratio comes from since it is a counting statistic rather than a chemical mechanism, and from a developmental and environmental theory.  In the composite language of these coupled basis theories the genetic sorts and the congruence $E$ become definable, and an idealized, deterministic Mendelian core becomes interpretable in the multiple lower theory in the sense of Definition~\ref{def:multiple_emergence}. The distinctions genetics draws are $\mathcal{L}_{\text{chem}}$-invisible not because chemistry is too coarse but because they live in a language chemistry does not contain, and multiple emergence seeks to add in the missing language. 

This is an important insight to emphasize. Genetics reduces to chemistry ontologically, in that every organism is a chemical system, while failing to be interpretable in chemistry epistemologically, even after chemistry is constrained, because its congruence and its characteristic relations are not constructible in the chemical language. Interpretability is regained only in an expanded, multiple lower theory. The statistical laws that give genetics its disciplinary content require, in addition, developing a measure over a larger space of models.

\section{Extensions of Emergence}
\label{sec:extensions}

Model theory is a very powerful tool for formalizing ``emergence'', much more broadly applicable
than any of the particular formalizations that have previously been investigated. Even so, there
are many instances of emergence that it cannot capture. In this section we present several
different extensions of model theory that allow us that capture far more 
instances of emergence.

An important property of several of these extensions is that they introduce probability measures
into the theories.
In the real world, outside of quantum mechanics, probability only arises due to uncertainty / ignorance
on the part of the scientist building the models. Probabilities are not part of the intrinsic, physical
processes. Accordingly, in the associated extensions of model theory below, the probability measures are all exogenous to the
theories, to be specified by the scientist using these extensions of model theory. 
An obvious variant would be to instead require in the usual way that the probability 
measures (and associated sigma algebras) can be defined as functions, in the model-theoretic sense, in
the relevant theories, and that those theories can prove that such functions are in fact proper
probability measures. In the interests of space, we do not consider such variants here.

\subsection{Probabilistic Emergence }
\label{sec:emerge_with_measure}
There are very many examples where one colloquially views a particular scientific discipline  that involves probability distributions as
``emerging'' from another discipline, that also involves probability distributions. In particular, very often the distributions defined in
the lower theory directly provide the distributions defined in the upper theory. We can capture this as follows. 

Suppose that a theory $T_A$ emerges from a theory $T_B$,
and so is interpretable into a theory $\widetilde{T_B}$ using translation $\Gamma = (\Delta, E, \{\varphi_s\}_{s \in \Sigma_A})$. Suppose as well that $\calL_A$ can be used to define a set $S_A$, and 
that each sentence in $T_A$ specifies that the variables in it must
be from $S_A$.\footnote{For example,
say that $S_A$ has $N$ elements, and that $\calL_A$ has $N$ constants,
$\{c_1, c_2, \ldots, c_N\}$. As shorthand, define the relation $R(x) = (x = c_1) \vee (x = c_2) \vee \ldots \vee (x = c_N)$. So $R$
defines the set $S_A$; 
$R(x)$ is true iff $x \in S_A$. Now suppose we are interested in 
a sentence that specifies that a relation $r(x, y) \in \calL_A$ is true 
if $x \le y$ for some total order $\le$ that is in the signature
$\calL_A$. To restrict such a
sentence to $S_A$, we would write it as
$R(x) \land R(y) \land [\le(x, y)] \rightarrow r(x, y)$.
}
Finally, suppose $T_A$ provides a sigma algebra over
$S_A$, $\mathscr{F}_{S_A}$. 

By~\cref{def:theory-interp}, there are formulas 
in $\calL_{\widetilde{T_B}}$ that define the equivalence
relation $\Delta / E$ and that say whether two 
associated equivalence classes are equal. In addition, 
it's required that the theory 
$\widetilde{T_B}$ proves that
those formulas define a proper equivalence relation.
Moreover, we can apply $\Gamma$ to the 
(extra-logical symbols in  the) formula in $\calL_A$ 
which defines $S_A$, to construct a bijection between $S_A$ and
a subset of $\Delta / E$ that we write as $\Gamma(S_A)$. Similarly, 
$\Gamma$ provides a sigma algebra isomorphism between 
$\mathscr{F}_{S_A}$ and a sigma algebra over
$\Gamma(S_A)$ that we write as $\Gamma(\mathscr{F}_{S_A})$.

Next, suppose we also have a sigma algebra $\mathscr{F}_{\widetilde{S_B}}$,
defined using formulas over $\calL_{\widetilde{T_B}}$,
and that $\mathscr{F}_{\widetilde{S_B}}$
is a refinement of the sigma algebra $\Gamma(\mathscr{F}_{S_A})$.
Then any probability measure P$_{\widetilde{T_B}}$ defined over 
the events defined by the sigma algebra $\calL_{\widetilde{T_B}}$
specifies a unique measure P$_{\Gamma(\mathscr{F}_{S_A})}$ over 
the events defined by $\Gamma(\mathscr{F}_{S_A})$.
That measure in turn specifies a unique measure over 
$\mathscr{F}_{S_A}$.


In terms of our framework, P$_{\widetilde{T_B}}$ 
is a probability measure over the states of any physical system 
that obeys the scientific formalism $T_B$ when there are auxiliary
constraints $\widetilde{T_B} \setminus T_B$. In turn,
$\Gamma$ provides a way to uniquely map that measure
to a probability measure over the states of certain physical
systems that obey the scientific formalism $T_A$. 

We summarize this with the following definition:

\begin{definition}(Probabilistic Emergence)
\label{def:emergence_measure}
Suppose that a theory $T_A$ emerges from a theory $T_B$,
using translation $\Gamma = (\Delta^B, E^B, \{\varphi_s\}_{s \in \Sigma_A})$. Suppose as well that
\begin{enumerate}
\item 
$\calL_A$ can be used to define a set $S_A$ and a sigma algebra over
that set, $\mathscr{F}_{S_A}$;
\item 
Each sentence in $T_A$ specifies that the variables in it must
be from $S_A$;
\item 
\label{item:with_measure_3}
$\calL_{\widetilde{T_B}}$ can be used
to define a sigma algebra $\mathscr{F}_{\widetilde{S_B}}$ over $\Delta^B$,
which is a refinement of the sigma algebra $\Gamma(\mathscr{F}_{S_A})$;
\item 
\label{item:with_measure_4}
We are given a probability measure P$_{\widetilde{T_B}}$
over the sigma algebra $\mathscr{F}_{\widetilde{S_B}}$.
\end{enumerate}
Finally, suppose as well that $\widetilde{T_B}$ proves both~\cref{item:with_measure_3} and~\cref{item:with_measure_4}.

In this situation P$_{\widetilde{T_B}}$
specifies a probability measure P$_{S_A}$ over
the sigma algebra $\mathscr{F}_{S_A}$, and we say
that $(T_A, S_A, \rm{P}_{S_A})$ 
\textbf{probabilistically emerges} 
from $(T_B, \Gamma(S_A), \rm{P}_{\widetilde{T_B}})$
(using translation $\Gamma = (\Delta, E, \{\varphi_s\}_{s \in \Sigma_A})$).
\end{definition}

Probabilistic emergence provides a way to 
have a probability measure in a lower theory $\widetilde{T_B}$ generate a probability measure in an an upper theory $T_A$ that emerges from $T_B$.
Loosely speaking, if we have a unary relation $R$ defined over $\calL_A$, then the distribution $P_{\Gamma(S_A)}$
in the lower theory, induced by the distribution $\rm{P}_{\widetilde{T_B}}$, tells us the probability that $R$ is true.
Moreover, each equivalencxe class
in $\Delta^B / E^B$ need not be an indivisible event of the sigma algebra $\mathscr{F}_{\widetilde{S_B}}$,
and there might be multiple events in such an equivalence class, which need not have the same probability.

As an illustration, the Boltzmann distribution over the
states of a finite system $C$ that is connected to an extremely large heat bath $D$
probabilistically emerges from the microcanonical distribution over the joint states of $C$ and $D$. We can use this, e.g., to evaluate
the probability that the energy of the system exceeds some given constant.

In general, the definition of the sigma algebra $\mathscr{F}_{\widetilde{S_B}}$ might use sentences that cannot
be constructed from the (appropriate quantifications of) the formulas that define $\Delta^B / E^B$. An
example where that is often the case is
ff there are microstates defined in both the lower theory and the upper theory, as discussed in~\cref{sec:interpretation_theories}, 
and we use those microstates to define the associated sigma algebras.

\subsection{Approximate Emergence}
\label{sec:approximate_emergence}

In many cases one formalism won't exactly arise from another one,
but only approximately do so, to within some accuracy. Typically such approximate 
emergence is combined with a probability measure over a sigma algebra defined using the
upper formalism's language. For example, often in gathering data we will have ``experimental''
or ``observational'' noise, that result in some ``error bars'' about the value of some observable
in the upper formalism. 

We can capture this scenario with the following extension of emergence:

\begin{definition}(Approximate emergence)
Suppose that $(T_A, S_A, \rm{P}_{S_A})$ 
{probabilistically emerges} from $(T_B, \Gamma(S_A), \rm{P}_{\widetilde{T_B}})$, and that it does so
using a translation that we write as $\Gamma = (\Delta^B, E^B, \{\varphi_s\}_{s \in \Sigma_A})$, along with
some implicit auxiliary constraints constructing $\widetilde{T_B}$ from $T_B$.
Suppose as well that we have a metric $d^B(., .)$ defined by appropriately combining
the formulas $\{\Gamma(\varphi_s)\}_{s \in \Sigma_A}$.  Suppose as well that ($\widetilde{T_B}$ proves that) $d^B(., .)$ is invariant over the 
equivalence classes $E^B / \Delta^B$. In this situation, $\Gamma^{-1}(d^B(.,.))$ is a well-defined formula in $\calL_A$, defining
a distance measure $d^A(., .)$.

Then we say that $(T_A, S_A, {\rm{P}_{S_A}}, d^A)$ \textbf{approximately emerges} 
from the associated tuple $(T_B, \Gamma(S_A), {\rm{P}_{\widetilde{T_B}}}, d^B)$
(using translation $\Gamma = (\Delta^B, E^B, \{\varphi_s\}_{s \in \Sigma_A})$).
\end{definition}

\noindent Note that the definition of approximate emergence implicitly assumes that there
is a set of cardinal numbers, e.g., $\N$ or $\Z$, definable in both $\widetilde{T_B}$ and $T_A$.
If this weren't the case, then it would not be possible for the distance measures to have numerical values.

Approximate emergence provides a way for us to use a provided metric
and probability measure over a lower theory to specify the
probability distribution over the distances
between any two elements (in the set $S_A$ defined)
in the upper theory.


\subsection{Stochastic Emergence}
\label{sec:stoch_emergence}

Often in cases of what is called ``emergence'' of one scientific formalism from another,
the laws themselves in the emergent formalism only hold
with a certain probability. This behavior is not possible using conventional model-theoretic
interpretation (assuming that the upper, emerging theory is consistent). 
Adding auxiliary constraints as in~\cref{def:emergence} doesn't salvage the situation.

To see how to capture the possibility of a law in an upper theory only holding
with a certain probability, note that in practice, in situations where we view a law in an 
emergent formalism only holding with some probability, what we really mean is that there
is a certain probability of its holding in a physical system which satisfies the lower formalism.

Given this,
suppose that $T_A$ emerges from $T_B$, but that $T_A$ is not
syntactically complete. So there are sentences $\phi$ with
signature $\calL_A$ that some models of $T_A$ disagree on. 
In particular, for a fixed translation $\Gamma$ from $T_A$ into
$T_B$, there may be two models $M$ and $N$ of $T_B$, where the associated models 
$\Gamma^M$ and $\Gamma^N$  disagree on $\phi$ (even though those two models 
are required to agree on all the sentences in $T_A$,
along with all logical consequences of those sentences). 
Given such a $\phi$, the emergence of $T_A$ from $T_B$ is not
faithful, and even if $\Gamma(\phi)$ is a logical consequence of $T_B$,
it doesn't follow that $\phi$ is satisfied in both $\Delta^M / E^M$ and $\Delta^N / E^N$. On the other hand, there might also be pairs of
$\calL_B$ structures $M', N'$ that both satisfy $T_B$ where
the structures $\Delta^{M'} / E^{M'}$ and $\Delta^{N'} / E^{N'}$
do agree on that same $\phi$. So whether $\phi$ is taken
to be true depends on what models of $T_B$ we are considering.

Now suppose that we have some measurable set $\mathscr{M}$ of models of
$T_B$ and an associated probability measure, $Pr$. If all four
models $M, N, M', N'$ are contained in $\mathscr{M}$, then 
there is an associated probability of whether the sentence $\phi$ in $\calL_A$ is true. This motivates the
following definition.

\begin{definition}(Stochastic Emergence)
\label{def:stochastic}
Suppose that a syntactically incomplete theory $T_A$ emerges from
a consistent theory $T_B$ using a translation
$\Gamma = (\Delta, E, \{\varphi_s\}_{s \in \Sigma_A})$. 
Suppose as well that $\mathscr{M}$
is a measurable set of models of $T_B$, and let $Pr_\mathscr{M}$ be a measure
over $\mathscr{M}$. 

Write $Pr_A$ for the probability measure over whether the
sentences in $\calL_A$ hold, induced by evaluating $Pr_\mathscr{M}$ 
over the structures $\{\Delta^M / E^M : M \in \mathscr{M}\}$ (which necessarily all
satisfy $T_A$).
We say that $(T_A, Pr_A)$ \textbf{stochastically
emerges} from $(T_B, Pr_\mathscr{M})$.
\end{definition}

Probabilistic emergence and stochastic emergence propagate probability in two complementary ways. Probabilistic emergence propagates a
 probability measure over the elements of a domain, from a lower theory $T_B$ to an upper theory $T_A$, and thereby assigns a probability to the event that a property, some relation in $\calL_A$, holds of a randomly chosen tuple of elements in the upper theory's
domain. Stochastic emergence instead propagates a measure over models of $T_B$, and thereby assigns a probability to the event that a law, some sentence of $\calL_A$, holds at all. In terms of our framework, the measure
$Pr_\mathscr{M}$ in stochastic emergence is a probability measure over a set of physical
systems, all governed by the scientific formalism $T_B$. It 
induces a probability measure over the possible logical
consequences of the axioms defining the formalism $T_A$ that
emerges from $T_B$.

A refined notion of stochastic emergence arises if we
restrict how uncertain the truth of any of the sentences in $\calL_A$
can be, under the measure $Pr_\mathscr{M}$:

\begin{definition}(Confident Emergence)
Suppose that $(T_A, Pr_A)$ {stochastically
emerges} from $(T_B, Pr_\mathscr{M})$. Suppose as well
that for some $\epsilon \in [0, 1/2)$,
$Pr_A(\phi) \not \in [\epsilon, 1 - \epsilon]$
for all sentences $\phi$ that can be formulated in $\calL_A$.
Then we say that $(T_A, Pr_A)$ \textbf{confidently 
emerges} from $(T_B, Pr_\mathscr{M})$
to within \textbf{precision} $\epsilon$.
\end{definition}

Confident emergence between theories may 
not be transitive, since errors may compound along chains of confident emergences. This provides
a formal reason why long sequences of emergence of one formalism from another degrade even when each link is nearly exact.

\subsection{Multiple Emergence}

It is not the case that (for example) biology emerges solely from chemistry.
It also involves Newtonian mechanics, e.g., to describe the tensile strength
of vertebrate bones. To capture this kind of phenomenon we need a way of having one formalism 
$\calB$ emerge from a set of multiple other formalisms,  $\{\calA_i\}$. That would mean 
that some formulas in $\calB$ can be interpreted in $\calA_1$, while some other formulas
in $\calB$ can be interpreted in $\calA_2$. In addition, often there are also coupling formulas,
that involve both $\calA_1$ and $\calA_2$. 
As an illustration
some of the formulas in biology
involve just chemistry, while others involve just Newtonian mechanics. However, there
are also formulas in biology that couple those formalisms, e.g., formulas 
coupling the Newtonian force on a particular bone to formulas governing
the ratio of osteoclasts to osteoblasts inside that bone.



We formalize this as follows:

\begin{definition} (Multiple emergence)
\label{def:multiple_emergence}
Let $\{\calA^i : i \in I\}$ be an indexed set of \textbf{basis} theories, with 
signatures $\{\calL_{i} : i \in I\}$. Let $\rho$ be a \textbf{coupling}
theory, with signature $\calL_\rho$.
Suppose we have an associated theory $\calB$ where:
\begin{enumerate}
\item $\calL_\calB = \calL_\rho \sqcup \coprod_{i \in I} \calL_{i}$;
\item $\calB$ is {multi-sorted}, in that the infinite set of
variables in the language of $\calB$ is partitioned,
into a collection of {sorts}, one sort for $\rho$ and one for each
theory $i \in I$. 
Write that set of variables of $\calB$ in terms of the sorts as
$X^{\calL_B} = X^\rho \cup \bigcup_{i \in I} X^i$;
\item For all $i \in I$, every relation or function
in $\calL_{i}$ only involves the variables in the 
associated sort, $X^i$;
\item For all pairs of basis theories $i, j \ne i$, there are functions or relations in
$\rho$ that involves both basis sorts $X^i$ and $X^j$.
\end{enumerate}
Then we say
that $\calB$ is \textbf{multiple interpretable} 
in the set of theories $\{\calA^i\}$ with
\textbf{coupling theory} $\rho$.
If some other theory $\calC$ is emergent from such a
$\calB$, then we say that $\calC$ is 
\textbf{multiple emergent} from the set of theories $\{\calA^i\}$ with
\textbf{coupling theory} $\rho$.
\end{definition}
\noindent Note that due to the third requirement in~\cref{def:multiple_emergence}, 
all functions in any one of the basis signatures are partial functions.

\begin{example}
\label{ex:electrodynamics}
Recall that Maxwell's equations specify the electric and magnetic fields throughout a region of space
based on the positions and velocities of all the particles in that region. In turn, the electromagnetic forces on any 
(possibly moving) charged particle are determined by the local electric and magnetic fields. 

Accordingly, electrodynamics -- the theory determining the dynamics of a Newtonian set of $N$
moving charged particles -- is multiple emergent from Newton's laws and Maxwell's equations. Concretely, $\calA^1$ is 
Newton's laws for a set of $N$ particles, and $\calA^2$ is Maxwell's equations. The associated coupling
theory $\rho$ provides the functions giving the Newtonian forces on the $N$ particles due to the electric and magnetic
fields. It also provides the function (just the identity function, in this trivial case) setting the charges
and velocities in Maxwell's equations to the values given by the associated quantities in Newton's laws.
Finally, the auxiliary constraints give the initial conditions of the $N$ particles.
\end{example}

\section{Relation with what's colloquially called ``emergence''}
\label{sec:applications}

In this section we show how our formalism captures some of the many phenomena
that are informally called ``emergence'' in the literature.

\subsection{Coarse-graining as equivalence classes in a lower theory}
\label{sec:coarse-graining}

Standard coarse-graining is a special case of the equivalence relation $E$ in Definition~\ref{def:struct-interp}. A partition of a state space $X$ into macroscopic bins $\{X_\alpha\}$ defines an equivalence relation, where $x \sim y$ iff $x$ and $y$ belong to the same bin. The macroscopic ``state'' is an equivalence class, and macroscopic observables are functions that are invariant under this equivalence, which is exactly the condition that model-theoretic interpretation requires of the translated symbols. This partition picture is a statement at the level of structures (Definition~\ref{def:struct-interp}), and a single microscopic state space is quotiented by $E$. The commutativity requirement introduced next is different in kind, and once it is demanded to hold uniformly across all models of the microscopic theory, it becomes a condition at the level of theories.

The canonical commutativity requirement considered in the literature
on what ``emergence''
means in the case of coarse-graining provides a further constraint. 
Given a microscopic dynamics $f: X \to X$ and a coarse-graining map $\pi: X \to Y$, the emergent dynamics $g: Y \to Y$ satisfies commutativity when $\pi \circ f = g \circ \pi$. In model-theoretic terms, the micro-level structure $\mathcal{B}$ has the dynamics $f$ in its signature, and the macro-level structure $\mathcal{A}$ has $g$ in its signature. Commutativity says that the interpretation map $\Gamma$ (which includes $\pi$ as the encoding from $A$ into equivalence classes of $B$) preserves the dynamical sentence, and the sentence ``$g$ maps state $\alpha$ to state $\beta$'' in $T_A$ translates to a sentence provable in $T_B$. Commutativity is thus a consequence of interpretability when the theories include dynamical axioms.

But commutativity fails for many natural coarse-grainings. Not every partition of a chaotic state space yields a Markov process on the partition elements. In model-theoretic language, the failure means that the proposed $\Gamma$ does not yield an interpretation, and the translated dynamical sentence is not provable from $T_B$. When it fails, one must either change the partition (change $E$) or add auxiliary constraints $C$ restricting to initial conditions for which the commutativity approximately holds. The latter is stochastic emergence.

\subsection{Symmetry breaking}
\label{sec:symmetry}

Symmetry breaking has two complementary readings in our framework. The collective modes of a large coupled system are an instance of multiple emergence, and the selection of one broken-symmetry ground state is an instance of stochastic emergence. Taking the collective modes first, recall that, as discussed briefly in~\cite{Anderson1972}, a very large system of spins
on a lattice can support propagation of phonons, i.e., waves in spin orientation that
move through the lattice. To model this
as a multiple emergent scenario, suppose we have a (thermodynamically) large
number of basis theories  $\calA^i$, all of which are identical. (So they all contain the 
same set of sentences in the same language.) These theories will represent the
individual spins in our system, with the sentences in each $\calA^i$ providing
those properties of spins that will be necessary for supporting the propagation
of phonons. 

Next, let the signature of all those theories each contain a special constant, which we call the ``position''. 
Finally, suppose our coupling theory $\rho$ assigns a separate location on the lattice to the position constant of each of 
the individual theories, by fixing them to lie on an actual lattice in Euclidean space. The coupling theory 
also provides the Hamiltonian function that couples each spin in the lattice to its neighbors, along with
the associated dynamic laws governing how each spin's state evolves as a function of the joint state of
its neighbors. Concretely, if the basis theories involve quantum mechanics, and the spins are not entangled,
then one of the sentences in $\rho$ will 
be Schrodinger's equation, operating over the product of the wave functions of the individual spins. 
If instead we are using a semi-classical approximation of the system of spins,
then the precise details of that approximation will provide the dynamics of the full system rather than
Schrodinger's equation..

Finally, let $T$ be a theory that is multiple emergent from those basis theories and that coupling theory,
with auxiliary constraints that specify the initial state of the spin system. $T$ contains sentences that
describe the propagation of phonons across large lattices. However, no set of sentences in the individual basis
provides a description of such a process. So there are dynamic properties of the full, large-scale system
that are completely absent from the individual, micro-scale subsystems. In particular, there are symmetry
relations that are obeyed by individual spins but not by a thermodynamically large set of
coupled spins. This is an example of Anderson had in mind when he wrote
that ``more is different'', and that ``broken symmetry'' arises as go from microscopic scale to
macroscopic condensed matter systems.

The selection among degenerate ground states is the second reading. A crystal obeys quantum electrodynamics, but QED does not determine which crystal structure will form, and the ground state breaks the rotational symmetry of the Hamiltonian. In the language of phase transitions, the Hamiltonian is parameterized by a control parameter (temperature $T$, magnetic field $h$), and below a critical value of the control parameter the system selects one of several degenerate ground states characterized by a nonzero order parameter (magnetization, crystal orientation).

Model theory captures part of this. An incomplete theory $T_B$ (such as the theory of a Hamiltonian at a given temperature) has multiple models that disagree on sentences not entailed by $T_B$. In the ferromagnet, $T_B$ at zero field below $T_c$ does not entail the sentence ``the magnetization is positive'' or its negation, and both are satisfied by some model of $T_B$. Adding the auxiliary constraint $C = \{\text{``magnetization} > 0\text{''}\}$ narrows the model class and makes the magnetized theory interpretable. The different choices of $C$ (positive versus negative magnetization) correspond to different broken-symmetry phases.

This is the paradigmatic case of stochastic emergence in the sense of Definition~\ref{def:stochastic}. The two broken-symmetry phases are not merely two admissible models, and in the absence of an external field they carry a symmetric measure, each occurring with probability one half. Writing $\mathscr{M}$ for the models of $T_B$ below $T_c$ and $Pr_{\mathscr{M}}$ for that measure, the magnetized theory holds in a measure-one set of models once one conditions on either phase, while the symmetric $\tfrac{1}{2}$--$\tfrac{1}{2}$ split over the two phases is itself the emergent-level datum that no single uniform constraint can recover. This is why symmetry breaking is stochastic emergence rather than strict interpretability, in that the split among models of $T_B$ is real, but it is governed by a measure rather than eliminated by an axiom.

What model theory, even equipped with a measure, does not supply is the dynamics of selection. Which ground state is realized depends on the system's history, on infinitesimal perturbations, and on boundary conditions. The framework classifies the possible outcomes and assigns them probabilities, but is silent on why a particular phase obtains on a given occasion. The control parameter $T$ and the order parameter $m$ are both expressible in $\mathcal{L}_B$ (they are definable functions on the domain), but the critical value $T_c$ at which the phase transition occurs is a property of the theory, and it is the value at which $T_B$ passes from completeness to incompleteness with respect to sentences about $m$. This locates phase transitions precisely in the model-theoretic architecture, as the points at which a parameterized family of theories loses completeness with respect to an order-parameter sentence.

\section{Computational complexity and the limits of interpretability}
\label{sec:complexity}


Model-theoretic interpretability of theories is an existence claim, asserting that there exists a special type of
translation $\Gamma$ that maps $T_A$ into $T_B$, obeying the conditions listed in~\cref{def:theory-interp}. 
It says nothing about the computational cost of constructing or evaluating $\Gamma$.

When the domains are countable and the translation $\Gamma$ is (Turing) computable, the width and
blow-up of the translation, defined in~\cref{sec:parsimony}, are measures of the compression that the emergent theory provides. 
We can define the usual computational complexity classes in terms of that blow-up. For example,
a sentence $\omega$ in $\mathcal{L}_A$ of length $n$ may translate to a sentence
$\Gamma(\omega)$ of length $n^k$ for some $k$ (polynomial) or $2^n$ (exponential). 
However, when dealing with finite theories, the blow-up $\widehat{\{\varphi_{s}\}}$ is finite, and so typically the number of symbols
in $\Gamma(\omega)$ is no more than $\widehat{\{\varphi_{s}\}}$ times the length $\omega$.
So in this case, the complexity is just linear.

A deeper question concerns the (Turing) decidability of the interpretation relation itself. 
Given two finitely axiomatized theories $T_A$ and $T_B$, can an algorithm determine 
whether $T_A \trianglelefteq T_B$? For broad classes of theories, the answer is no~\cite{Visser2006}. 
Even when an interpretation exists, there may be no algorithm for discovering it. We cannot in 
general compute the lattice of interpretability degrees, and can only explore it piecemeal.

Sometimes an upper, ``macroscopic'' theory that involves dynamics  emerge
from a lower, ``microscopic''  that also involves dynamics. Moreover, sometimes there
are properties of the macroscopic dynamics that can only be determined by simulating the full microscopic dynamics step by step. 
In our framework, this corresponds to an interpretation $\Gamma$ that exists but whose evaluation on a given 
model requires time proportional to the full simulation. The interpretation is logically valid but computationally 
useless, and knowing that $T_A \trianglelefteq T_B$ does not help you predict the behavior of a particular model 
of $T_B$ any faster than brute simulation.

This connects to a question raised by the programming analogy. A compiler maps a high-level program 
into machine code, and the map runs in polynomial time. But the execution of the resulting machine code, 
on a particular structure (a particular chip with particular inputs), may take exponential time. Knowing that $T_A \trianglelefteq T_B$ tells you only that a translation exists and is cheap to write down, and it says nothing about how expensive that translation is to evaluate on a given model. Emergence is the cheap, theory-level fact, and simulation is the expensive, model-level one. Model theory certifies the former and is silent on the latter.

The interested reader is directed to the literature on finite model theory, in which 
issues concerning the computational complexity of interpretations plays a central role~\cite{libkin2004elements,mccolm1996ebbinghaus}.

\section{Model Theory in The Wild}
\label{sec:engineered-natural}

The programming example of Section~\ref{sec:programming} is misleading as an exemplar of scientific emergence, because the compiler was designed, and its interpretation map $\Gamma$, its equivalence relation $E$, and its auxiliary constraints $C$ were all specified in advance rather than discovered. Most scientific emergence is not like this. No one designed the mapping between neural activity and mental representation, or between individual behavior and collective dynamics. If those maps can be formalized as interpretations at all, they are products of evolution, development, or self-organization. We examine two such cases, brain and mind, and individuals and collectives, to see where the formalism marks the difference between evolution and engineering.

\subsection{Brain and mind}

Consider the claim that a theory of mental states $T_{\text{mind}}$ emerges from a theory of neural activity $T_{\text{neural}}$. In the model-theoretic framework, this requires a translation $\Gamma$ that constructs the sorts, functions, and relations of $T_{\text{mind}}$ from formulas in $\mathcal{L}_{\text{neural}}$.

In the compiler case, the signatures of both $T_H$ and $T_M$ are fixed before the interpretation begins, along with 
the associated theories. We know what ``integer addition'' means in $\mathcal{L}_H$, i.e., we know the sentences that describe how the symbol
for integer addition combines with the symbols for other concepts in the high level language. Similarly, we know what ``word-level ADD'' means in $\mathcal{L}_M$. The compiler translates the one into the other. 

For brain and mind the associated signatures aren't settled. In particular, though it's not phrased this way,
the signature of $T_{\text{mind}}$ (the non-logical symbols for belief, desire, intention, percept, affect) is a matter of ongoing dispute among cognitive scientists, philosophers, and psychologists. Different theoretical traditions assign different arities to these symbols, disagree about whether they are function symbols or relation symbols, and in some cases deny their existence altogether. The signature of $T_{\text{neural}}$ is comparatively better defined (firing rates, synaptic weights, receptor densities, network topology) but is itself a moving target as new measurement technologies (calcium imaging, Neuropixels, connectomics) reveal structure at finer resolutions.

When neither signature is fixed, the interpretation problem changes character. One is not translating a given vocabulary into another given vocabulary, but simultaneously discovering both vocabularies and the mapping between them. In engineered systems, signatures precede interpretation. In natural systems, the act of interpretation participates in the construction of the signatures. (Never mind the associated sentences.) This is why in previous work we have described Reality as an Oroborous, a snake that eats its own tail. In order to make foundational claims in most complex systems, one needs to borrow from a "higher" or composite measurement, thereby shifting what is foundational to another levels of organization. 

In the compiler, the equivalence relation $E$ in the lower theory
is easily stated as two bit-patterns either encode the same integer or they do not. There is no ambiguity and no context-dependence. 

Neural equivalence classes are very different.
What counts as ``the same mental state'' across different patterns of neural activity is not a binary relation that can be read off the neural domain. The neuroscience of multi-voxel pattern analysis and neural decoding all involve imposing equivalence structures on high-dimensional activation spaces, but these structures are observer-dependent (they require a choice of similarity metric), context-dependent (the same neuron may participate in different equivalence classes depending on the task), and they are approximate (the boundaries of the classes are graded). In the model-theoretic formalism, $E$ must be a definable equivalence relation on $B^k$ for some $k$. If the ``correct'' equivalence relation is not definable in $\mathcal{L}_{\text{neural}}$ (because it depends on the organism's history, or on environmental context that lies outside the neural domain) then the interpretation cannot be carried out within the neural theory alone. One needs an expanded lower theory that includes the environment, the developmental history, or both. The auxiliary constraints $C$ are doing much heavier work than in the engineered case, and it is not clear that they can be finitely axiomatized.

In the engineered case, the interpretation $\Gamma$ is an explicit, surveyable object. One can inspect the compiler's source code, verify its correctness, trace its action on any particular input. For brain and mind, no such object exists. The ``compiler'' is the product of roughly 600 million years of nervous system evolution and roughly 20 years of individual neural development. It is not written in any inspectable language. Its existence is conjectured on the basis of correlational evidence (neural imaging, lesion studies, optogenetics) and theoretical argument (functionalism, computational theories of mind), but the actual map $\Gamma$, if it exists, has never been exhibited.

Model theory is indifferent to this gap. Interpretability is an existence claim, asserting that there exists a translation satisfying certain properties. Whether anyone has constructed the translation, or whether it is constructible in principle by a finite agent, is not a question the formalism addresses.  The claim that $T_{\text{mind}} \trianglelefteq T_{\text{neural}} \cup C$ may be true yet permanently beyond our ability to verify, not because of computational cost (the worry of Section~\ref{sec:complexity}) but because the interpretation map is not the sort of thing that admits explicit construction. It was assembled by a historical process, not a logical one. 

The extensions to model theory introduced above clarify much of what makes the brain-mind relation harder than the compiler case. The compiler case illustrates one kind of multiple realizability directly, and the same interpretation $\Gamma$, proved from $T_M$, constructs $T_H$ inside every chip. Different chips are different structures $\mathcal{M} \models T_M$, so the variability is across models of a fixed $\Gamma$. The brain case involves a different kind of variability. Different people may encode the same concept through quite different patterns of neural activity, so there is no single $\Gamma$ that holds across individuals, and there is instead a family $\{\Gamma_i\}$ of candidate interpretations, one per individual, with the variability across the interpretations themselves rather than across models of one interpretation. This is exactly where the parallel to the compiler breaks down, and it is the confound to avoid, since the object that is stable across individuals is not $\Gamma$ but the emergent theory $T_{\text{mind}}$ that each $\Gamma_i$ induces.

With the confound removed, the extensions slot in cleanly, and this is the constructive payoff. Individual variation in $\Gamma$ is captured by stochastic emergence, and placing a measure over the family $\{\Gamma_i\}$, equivalently over individuals (models of $T_{\text{neural}}$), makes the emergent mental laws hold with probabilities, and ``the same concept, different neural codes'' becomes a distribution over interpretations. The graded, similarity-based character of mental states, evident in the multi-voxel and decoding analyses above, is captured by approximate emergence, with the decoding metric supplying the metric $d^A$ on emergent states. And the dependence of the ``correct'' equivalence relation $E$ on environment and developmental history, rather than being a matter of auxiliary constraints that cannot be finitely axiomatized, is captured by multiple emergence, and $T_{\text{mind}}$ emerges jointly from $T_{\text{neural}}$ and a coupled environmental or developmental theory (Definition~\ref{def:multiple_emergence}), with the missing structure living in a second basis theory rather than in an unbounded axiom list.

What none of these extensions supplies is a constructed $\Gamma$. Model theory and its probabilistic extensions certify the existence and structure of the interpretation, and they never exhibit it by a finite agent. The claim that $T_{\text{mind}}$ emerges from $T_{\text{neural}}$ may therefore be true yet permanently beyond our ability to verify, because the interpretation was assembled by a historical process rather than a logical one. 

\subsection{Individuals and collectives}

Consider now the claim that a theory of collective behavior $T_{\text{coll}}$ emerges from a theory of individual behavior $T_{\text{ind}}$.

In the compiler case, the ``parts'' of the machine (registers, ALU, memory cells) interact according to a fixed architecture specified by the hardware design. The interaction structure is part of the signature of $T_M$, and there are function symbols for reading from and writing to specified addresses, with the connectivity hard-wired.

In a collective (a colony of ants, a flock of birds, a market of traders, a population of neurons) the interaction structure is itself a variable. Who interacts with whom, when, and through what channel is not fixed by the theory of individual behavior, and it is an emergent property of the dynamics. The network of interactions is not in $\Sigma_{\text{ind}}$, and it arises from the interplay of individual strategies, spatial proximity, signalling, and environmental constraints. In model-theoretic terms, the interaction topology belongs not to $T_{\text{ind}}$ but to the auxiliary constraints $C$ or perhaps to an expanded signature.

This creates a bootstrapping problem. To interpret $T_{\text{coll}}$ in $T_{\text{ind}} \cup C$, we need $C$ to specify the interaction structure. But the interaction structure is itself a collective-level property, part of what $T_{\text{coll}}$ describes. The auxiliary constraints that enable emergence are drawn from the vocabulary of the emergent theory, not just from the vocabulary of the lower theory. This circularity has no obvious counterpart in the engineered case, where the hardware architecture is given independently of the software that runs on it.

What counts as ``the same collective state'' across different configurations of individuals? In a bird flock, two configurations that differ in the positions of individual birds but share the same velocity field, density profile, and polarization may be treated as equivalent. In a market, two microeconomic states that produce the same price, volume, and volatility may be regarded as the same macroeconomic state. These equivalence relations are not arbitrary, and they are selected by the collective dynamics itself. The macroscopic observables that define the equivalence classes are precisely those that are approximately preserved under the microscopic dynamics, and these are the slowly varying modes, the conserved quantities, and the order parameters.

In model-theoretic language, the equivalence relation $E$ is not given in advance (as in the compiler) or discovered by correlation (as in the brain-mind case), and it is selected by the dynamics from among the possible partitions of the microstates space. The correct $E$ is the one for which the interpretation $\Gamma$ yields a well-defined emergent dynamics satisfying commutativity (Section~\ref{sec:coarse-graining}). This means that the equivalence relation, the translation map, and the emergent theory are jointly determined, and they are not independent components to be specified one at a time, as they are in the engineered case. The interpretability problem in collective systems is a fixed-point problem, and one must find $E$, $\Gamma$, and $T_{\text{coll}}$ simultaneously such that the interpretation holds.

In the compiler case, $T_M$ is fixed, and the instruction set architecture does not change depending on what program is being compiled. In collective systems, the effective lower theory can shift as the collective organizes. Individual ants in a colony modify their behavioral repertoire in response to colony-level states (task allocation, alarm signals, food availability). Individual neurons modify their synaptic weights in response to network-level activity. The lower theory $T_{\text{ind}}$ is not a fixed set of axioms, and it is itself modulated by the emergent level. In model-theoretic terms, the axioms of $T_{\text{ind}}$ are functions of the model of $T_{\text{coll}}$ in which the individual is embedded.

\section{A Concluding Typology of Emergence}
\label{sec:concluding_typology}

Recapping, in our framework, we identify the mathematically stated laws of a particular scientific formalism
with a particular first order theory. A canonical example is Newton's laws, by themselves.
Note that this formalism does not specify any forces, masses, phase space positions, etc. 

We then identify any particular physical system that obeys such a set of laws with a structure 
that satisfies that theory. That structure is where the forces, masses, and phase space positions
are specified.  A system of gravitating point masses in vacuum evolving
according to Newton's laws is an example.

If an upper formalism $T_A$ is interpretable in a lower formalism $T_B$ but not vice-versa,
then we identify $T_A$ as emerging from $T_B$. More generally though, we also
view $T_A$ as emerging from $T_B$ if $T_A$ is not interpretable in $T_B$, but it is
interpretable in $T_B$ together with a finite set of auxiliary constraints, sentences of $\calL_B$ that we add to $T_B$. This is the case for
example when we view statistical mechanics of ideal gases as emerging from Newton's laws,
and we need to add special auxiliary constraints concerning the probability measure
over the possible joint phase space positions of the particles in the gas in order to
get statistical mechanics.  

In addition, we use the term scientific discipline to refer to a set of formalisms that are
mutually interpretable. Classical dynamics is an example, which comprises the formalisms of
Newton's laws, Lagrangian dynamics, Hamilton's equations, etc. Since mutual interpretability
is transitive, all formalisms in a given discipline emerge from one another.

The several phenomena that go under the name ``emergence'' can now be organized through a model-theoretic typology. A typology does not resolve the ambiguities of emergence, but it makes explicit the role of auxiliary constraints and of a measure over the models of the lower theory, and it shows that otherwise disparate accounts are manifestations of the same logical relation.

Strict interpretability ($T_A \trianglelefteq T_B$, no auxiliary constraints) covers cases where the emergent theory is fully constructible from the lower theory. The compiler example belongs here, as does the interpretation of the rationals inside the integers.

Emergence ($T_A \trianglelefteq T_B \cup C$, with auxiliary constraints) covers the typical case in physics, and includes thermodynamics from statistical mechanics, hydrodynamics from molecular dynamics, and effective field theories from more fundamental ones. The constraints $C$ encode typicality, boundary conditions, or scale separation. These all add additional parameters to theories leading to significant departures from parsimony. 

Stochastic emergence covers cases where the emergent regularity holds with high probability under a natural measure. The ideal gas law holds for typical configurations in the microcanonical ensemble, and atypical configurations exist but have measure zero in the thermodynamic limit. Symmetry breaking is the other paradigm case, and below the critical temperature the broken-symmetry phases carry a measure, with the emergent theory holding in a measure-one set of models once a phase is fixed. This is different from the simpler idea of a probabilistic emergence function where the laws are fixed but the outcomes are distributions. 

Multiple emergence covers the special sciences, which emerge from conjunctions of lower theories (biology from chemistry, mechanics, and information theory jointly) rather than from any single base. The relationship between biology and physics is of this kind. Organisms satisfy biological laws, but not all physical systems are organisms, and no finite axiom set picks out the biological ones, so biology emerges from physics jointly with organizational and selective constraints rather than from physics alone. This is also described as exaptation in the case where constraints can be apportioned to either a system (e.g. traits in an organisms) and to the environment (e.g. phases of matter). Hence  a lung (providing a physiological interface for gas exchange) is an exaptation from a swim bladder (conferring neutral buoyancy in a fluid). We might also say that the mathematical abilities of current Large Language Models are exaptations from their original linguistic competence as natural language models. The change in environment or context in the case of an LLM is a change in syntax. 

Expressive dependence and foundational dependence are not emergence at all, but structural relations between theories (expansion, inclusion, identity) that model theory also formalizes.

Realizational dependence is emergence applied to the mind-body problem, where the question of whether the cognitive sorts are definable in neural vocabulary remains open.

Epistemological dependence is a constraint on the admissible class of theories, not an interpretability relation between theories, and may require resources beyond first-order model theory to formalize.

\section*{Acknowledgments}

DHW thanks for the Santa Fe Institute for support. DCK was supported by John Templeton Foundation Award 63797 on The Emergence of Agentic Intelligence. We would like to thank Sean Carroll and Joel Hamkins for helpful discussion.

\newpage

\bibliographystyle{plain}
\bibliography{references,references_additions}

\end{document}